\documentclass[aap,MSNbibl,seceqn,citesort,dvips]{arximspdf}
\usepackage{mathbh}
\usepackage{graphicx}

\doi{10.1214/11-AAP773}
\volume{22}
\issue{2}
\pubyear{2012}
\firstpage{541}
\lastpage{575}

\makeatletter

\newtheorem{Thm}{Theorem}
\newtheorem{prop}{Proposition}[section]
\newtheorem{lemme}[prop]{Lemma}
\newtheorem{corol}[prop]{Corollary}

\newproclaim{hyp}[prop]{Assumption}
\newproclaim{Rque}{Remark}
\newproclaim{ex}{Example}

\newcommand{\Co}{\mathcal{C}}
\newcommand{\supp}{\operatorname{supp}}
\newcommand{\rmd}{\mathrm{d}}

\newcommand{\D}{\mathcal{D}}
\newcommand{\B}{\mathcal{B}}
\newcommand{\C}{\mathcal{C}}
\newcommand{\U}{\mathcal{U}}
\newcommand{\N}{\mathbb{N}}
\newcommand{\PP}{\mathbb{P}}
\newcommand{\R}{\mathbb{R}}
\newcommand{\E}{\mathbb{E}}

\newcommand{\I}{\mbox{\textup{\textsc{i}}}}
\newcommand{\II}{\mbox{\fontsize{6.6pt}{6.6pt}\selectfont{\textup{I}}}}
\newcommand{\scS}{\mbox{\textup{\textsc{s}}}}
\newcommand{\sscS}{\mbox{\fontsize{6.6pt}{6.6pt}\selectfont{\textup{S}}}}
\newcommand{\Rr}{\mbox{\textup{\textsc{r}}}}
\newcommand{\sRr}{\mbox{\fontsize{6.6pt}{6.6pt}\selectfont{\textup{R}}}}
\newcommand{\cS}{\mbox{\textup{\textsc{s}}}}
\newcommand{\cI}{\mbox{\textup{\textsc{i}}}}
\newcommand{\cRr}{\mbox{\textup{\textsc{r}}}}
\newcommand{\scIS}{\mbox{\fontsize{6.6pt}{6.6pt}\selectfont{\textup{IS}}}}
\newcommand{\cRSr}{\mbox{\fontsize{6.6pt}{6.6pt}\selectfont{\textup{RS}}}}

\newcommand{\bpI}{\bar{p}^{{\II}}}
\newcommand{\bpS}{\bar{p}^{\sscS}}
\newcommand{\bpR}{\bar{p}^{\sRr}}
\newcommand{\muS}{\mu^{\sscS}}
\newcommand{\muSI}{\mu^{\scIS}}
\newcommand{\muSR}{\mu^{\cRSr}}
\newcommand{\muSn}{\mu^{(n),\sscS}}
\newcommand{\muSIn}{\mu^{(n),\scIS}}
\newcommand{\muSRn}{\mu^{(n),\cRSr}}
\newcommand{\muISn}{\mu^{(n),\scIS}}
\newcommand{\muRSn}{\mu^{(n),\cRSr}}
\newcommand{\muSnn}{\mu^{n,\sscS}}
\newcommand{\muSInn}{\mu^{n,\scIS}}
\newcommand{\muISnn}{\mu^{n,\scIS}}
\newcommand{\muRSnn}{\mu^{n,\cRSr}}
\newcommand{\bmuS}{\bar{\mu}^{\sscS}}
\newcommand{\bmuSI}{\bar{\mu}^{\scIS}}
\newcommand{\bmuSR}{\bar{\mu}^{\cRSr}}
\newcommand{\bmuIS}{\bar{\mu}^{\scIS}}
\newcommand{\bmuRS}{\bar{\mu}^{\cRSr}}
\newcommand{\NIS}{N^{\II\sscS}}
\newcommand{\NS}{N^{\sscS}}
\newcommand{\NRS}{N^{\sRr\sscS}}
\newcommand{\bNIS}{\bar N^{\II\sscS}}
\newcommand{\bNS}{\bar N^{\sscS}}
\newcommand{\bNRS}{\bar N^{\sRr\sscS}}
\newcommand{\NISn}{N^{(n),\II\sscS}}
\newcommand{\NSn}{N^{(n),\sscS}}
\newcommand{\NRSn}{N^{(n),\sRr\sscS}}
\newcommand{\NISnn}{N^{n,\II\sscS}}
\newcommand{\NSnn}{N^{n,\sscS}}
\newcommand{\NRSnn}{N^{n,\sRr\sscS}}
\newcommand{\In}{I^{(n)}}
\newcommand{\Sn}{S^{(n)}}
\newcommand{\Rn}{R^{(n)}}
\newcommand{\M}{\mathcal{M}}
\newcommand{\car}{\mathbf{1}}

\newcommand{\ind}{{\mathbh1}}

\makeatother

\begin{document}
\begin{frontmatter}

\title{Large graph limit for an SIR process in random network with heterogeneous connectivity}
\runtitle{Large graph limit for an SIR process}

\begin{aug}
\author[A]{\fnms{Laurent} \snm{Decreusefond}\ead[label=e1]{laurent.decreusefond@telecom-paristech.fr}\ead[label=u1,url]{http://perso.telecom-paristech.fr/\textasciitilde decreuse/}},
\author[B]{\fnms{Jean-St\'{e}phane} \snm{Dhersin}\thanksref{t1}\ead[label=e2]{dhersin@math.univ-paris13.fr}\ead[label=u2,url]{http://www.math.univ-paris13.fr/\textasciitilde dhersin/}},\\
\author[C]{\fnms{Pascal} \snm{Moyal}\ead[label=e3]{pascal.moyal@utc.fr}\ead[label=u3,url]{http://www.lmac.utc.fr/\textasciitilde moyalpas/}} and
\author[D]{\fnms{Viet Chi} \snm{Tran}\corref{}\thanksref{t1}\ead[label=e4]{chi.tran@math.univ-lille1.fr}\ead[label=u4,url]{http://labomath.univ-lille1.fr/\textasciitilde tran}}
\runauthor{Decreusefond, Dhersin, Moyal and Tran}
\affiliation{Telecom Paristech, Universit\'{e} Paris 13,
Universit\'{e} de Technologie de Compi\`{e}gne, and
Universit\'{e} des Sciences et Technologies Lille 1 and Ecole Polytechnique}
\address[A]{L. Decreusefond\\
Institut Telecom\\
Telecom Paristech, CNRS LTCI\\
23 av. d'Italie \\
75013 Paris\\
France\\
\printead{e1}\\
\printead{u1}}
\address[B]{J.-S. Dhersin\\
D\'{e}partement de Math\'{e}matiques\\
Institut Galil\'{e}e, LAGA, UMR 7539\\
Universit\'{e} Paris 13\\
99 av. J.-B. Cl\'{e}ment\\
93430 Villetaneuse\\
France\\
\printead{e2}\\
\printead{u2}}
\address[C]{P. Moyal\\
Laboratoire de Maths Appliqu\'{e}es\\
\quad de Compi\`{e}gne\\
Universit\'{e} de Technologie de Compi\`{e}gne\\
D\'{e}partement G\'{e}nie Informatique \\
Centre de Recherches de Royallieu \\
BP 20 529 \\
60205 Compi\`{e}gne Cedex\\
France\\
\printead{e3}\\
\printead{u3}}
\address[D]{V. C. Tran\\
Laboratoire Paul Painlev\'{e}\\
UMR CNRS 8524\\
UFR de Math\'{e}matiques\\
Universit\'{e} des Sciences\\
\quad et Technologies Lille 1\\
59655 Villeneuve d'Ascq Cedex\\
France\\
\printead{e4}\\
\printead{u4}\\
and\\
Centre de Math\'{e}matiques Appliqu\'{e}es\hspace*{9.6pt}\\
UMR 7641\\
Ecole Polytechnique\\
Route de Saclay \\
91128 Palaiseau Cedex\\
France}
\end{aug}

\thankstext{t1}{Supported by ANR Viroscopy (ANR-08-SYSC-016-03)
and MANEGE (ANR-09-BLAN-0215), and from the ``Chair Mod\'{e}lisation
Math\'{e}matique de la Biodiversit\'{e} of Veolia Environnement-Ecole
Polytechnique-Museum National d'Histoire Naturelle-Fondation~X.''}

\received{\smonth{10} \syear{2010}}
\revised{\smonth{2} \syear{2011}}

%
\begin{abstract}
We consider an SIR epidemic model propagating on a configuration model
network, where the degree distribution of the vertices is given and
where the edges are randomly matched. The evolution of the epidemic is
summed up into three measure-valued equations that describe the degrees
of the susceptible individuals and the number of edges from an
infectious or removed individual to the set of susceptibles. These
three degree distributions are sufficient to describe the course of the
disease. The limit in large population is investigated. As a corollary,
this provides a rigorous proof of the equations obtained by Volz
[\textit{Mathematical Biology} \textbf{56} (2008) 293--310].
\end{abstract}

%

\end{frontmatter}

\section{Introduction and notation}

In this work, we investigate an epidemic spreading on a random graph
with fixed degree distribution and evolving according to an SIR
model as follows. Every individual not yet infected is assumed to be
susceptible. Infected individuals stay infected during random
exponential times with mean $1/{\beta}$ during which they infect each
of their susceptible neighbors with rate $r$. At the end of the
infectious period, the individual becomes removed and is no longer
susceptible to the disease. Contrary to the classical
mixing compartmental SIR epidemic models (e.g.,
\cite{kermackmckendrick,bartlett} and see \cite{andersonbritton}, Chapter
2, for a
presentation), heterogeneity in the number of contacts makes it
difficult to describe the dynamical behavior of the epidemic. Mean
field approximations (e.g.,
\cite
{pastorsatorrasvespignani,barthelemybarratpastorsatorrasvespignani,durrett})
or large
population approximations (e.g., \cite{ballneal}, see also
equation (3) of
\cite{andersson} in discrete time) provide a set of denumerable
equations to describe our system. We are here inspired by the paper of
Volz \cite{volz}, who proposes a low-dimensional system of five
differential equations for the dynamics of an SIR model on
a~configuration model (CM) graph \cite{bollobas2001,molloyreed}. We refer
to Volz's article for a~bibliography about SIR models on graphs (see
also Newman \cite{newman,newmanSIAM}, Durrett~\cite{durrett} or
Barth\'{e}lemy et al.
\cite{barthelemybarratpastorsatorrasvespignani}). Starting from a
random model in finite population, Volz derives deterministic
equations by increasing the size of the network, following in this
respect works of Newman, for instance, \cite{newmanSIAM}. The
convergence of the continuous-time stochastic SIR model to
its deterministic limit for large graphs was, however, not
proved. In this paper, we prove the convergence that was left open by
Volz. To achieve this, we provide a~rigorous individual-based
description of the epidemic on a random graph. Three degree
distributions are sufficient to
describe the epidemic dynamics. We describe these distributions by
equations in the space of measures on the set of nonnegative integers,
of which Volz's equations
are a by-product. Starting with a node-centered description, we show
that the individual dimension is lost in the large graph limit. Our
construction heavily relies on the choice of a CM for
the graph underlying the epidemic, which was also made in \cite
{volz}.

The size $N$ of the population is fixed. The individuals are related
through a random network and are represented by the \textit{vertices}
of an undirected graph. Between two neighbors, we place an \textit
{edge}. The graph is nonoriented and an edge between $x$ and $y$ can be
seen as two directed edges, one from~$x$ to $y$ and the other from $y$
to $x$. If we consider an edge as emanating from the vertex $x$ and
directed to the vertex $y$, we call $x$ the \textit{ego} of the edge
and $y$ the \textit{alter}. The number of neighbors of a given
individual is the \textit{degree} of the associated vertex. The degree
of $x$ is denoted $d_x$. It varies from an individual to another one.
The CM developed in Section \ref{secconfigurations-graph} is a random
graph where individuals' degrees are independent random variables with
same distribution $(p_k)_{k\in\N}$. Edges are paired at random. As a
consequence, for a given edge, alter has the size-biased degree
distribution: the probability that her degree is $k$ is $k p_k/\sum
_{\ell\in\N}\ell p_\ell$.

The population is partitioned into the classes of susceptible,
infectious or removed individuals. At time $t$, we denote by $\scS_t$,
$\I_t$ and $\Rr_t$ the set of susceptible, infectious and removed
nodes. We denote by $S_t$, $I_t$ and $R_t$ the sizes of these classes
at time $t$. With a slight abuse, we will say that a~susceptible
individual is of type $\scS$ (resp., of type $\I$ or $\Rr$) and
that an edge linking an infectious ego and susceptible alter is of type
$\cI\cS$ (resp., $\cRr\cS$, $\cI\cI$ or~$\cI\cRr$). For
$x\in\cI$ (resp., $\cRr$), $d_x(\cS)$ represents the number
of edges with $x$ as ego and susceptible alter. The numbers of edges
with susceptible ego (resp., of edges of types $\cI\cS$ and $\cRr\cS
$) are denoted by $N^{\sscS}_t$ (resp., $\NIS_t$ and $\NRS_t$).

A possible way to rigorously describe the epidemics' evolution is given
in Section \ref{sectiondynamics}. We consider the subgraph of
infectious and removed individuals with their degrees. Upon infection,
the infectious ego chooses the edge of a~susceptible alter at random.
Hence, the latter individual is chosen proportionally to her degree.
When she is connected, she uncovers the edges to neighbors that were
already in the subgraph.

We denote by $\N$ the set of nonnegative integers and by $\N^*=\N
\setminus\{0\}$. The space of real bounded functions on $\N$ is
denoted by $\B_b(\N)$. For any $f \in\B_b(\N)$, set $\| f \|
_{\infty}$ the supremum of
$f$ on $\N$. For all such $f$ and $y\in\N$, we denote by $\tau_yf$
the function $x \mapsto f(x-y)$. For all $n\in\N$, $\chi^n$ is the
function $x \mapsto x^n$, and in particular, $\chi\equiv\chi^1$ is
the identity function and
$\mathbf1 \equiv\chi^0$ is the function constantly equal to 1.

We denote by $\M_F(\N)$ the set of finite measures on $\N$, embedded
with the topology of weak convergence. For all $\mu\in\M_F(\N)$ and
$f\in\B_b(\N)$, we write
\[
\langle\mu,f \rangle=\sum_{k\in\N}f(k) \mu(\{k\}).
\]
With some abuse of notation, for all $\mu\in\M_F(\N)$ and $k \in\N
$, we denote $\mu(k)=\mu(\{k\})$. 
For $x\in\N$, we write $\delta_x$ for the Dirac measure at point $x$.
Note, that some additional notation is provided in the \hyperref
[appendiceA]{Appendix}, together with several topological
results, that will be used in the sequel.

The plan of the paper and the main results are described below. In
Section~\ref{sectionmodel}, we describe the mechanisms underlying the
propagation of the epidemic on the CM graph. To describe the course of
the epidemic, rather than the sizes $S_t$, $I_t$ and~$R_t$, we consider
three degree distributions given as point measures of $\M_F(\N)$, for
$t\geq0$:
%
%
\begin{eqnarray}\label{muS}
\muS_t&=&\sum_{x\in\sscS_t}\delta_{d_x}
,\qquad \muSI_t
=\sum_{x\in
\II_t}\delta_{d_x(\sscS_t)},\nonumber\\[-8pt]\\[-8pt]
\mu^{\sRr\sscS}_t&=&\sum_{x\in
\sRr_t}\delta_{d_x(\sscS_t)}.\nonumber
\end{eqnarray}
Notice that the measures $\muS_t/S_t$, $\muSI_t/I_t$ and $\muSR_t/R_t$
are probability measures that correspond to the usual (probability)
degree distribution. The degree distribution $\muS_t$ of susceptible
individuals is needed to describe the degrees of the new infected
individuals. The measure $\muSI_t$ provides information on the number
of edges from $\cI_t$ to $\cS_t$, through which the disease can
propagate. Similarly, the measure $\muSR_t$ is used to describe the
evolution of the set of edges linking $\cS_t$ to $\cRr_t$.
We can see that $N^{\sscS}_t=\langle\muS_t,\chi\rangle$ and
$S_t=\langle\muS_t,\car\rangle$ (and, accordingly, for $\NIS_t$,
$\NRS$, $I_t$ and $R_t$).

In Section \ref{sectionlargegraph}, we study the large graph limit
obtained when the number of vertices tends to infinity, the
degree distribution being unchanged. The degree distributions
mentioned above can then be approximated, after proper scaling, by the solution
$(\bmuS_t,\bmuSI_t,\bmuSR_t)_{t\geq0}$ of
the system of deterministic measure-valued equations (\ref
{limitereseauinfiniS})--(\ref{limitereseauinfiniSR}) with initial
conditions $\bmuS_0$, $\bmuSI_0$\vspace*{2pt} and~$\bmuRS_0$.

For all $t\ge0$, we denote
by $\bNS_t=\langle\bmuS_t,\chi\rangle$ (resp., $\bNIS_t=\langle
\bmuSI_t,\chi\rangle$ and $\bNRS_t=\langle\bmuRS_t,\chi\rangle
$) the continuous number of edges with ego in $\scS$ (resp., $\I\scS$
edges, $\Rr\scS$ edges). Following Volz \cite{volz}, pertinent
quantities are the proportions
$\bpI_t=\bNIS_t/\bNS_t$ [resp., $\bpR_t=\bNRS_t/\bNS_t$ and $\bpS
_t=(\bNS_t-\bNIS_t-\bNRS_t)/\bNS_t$]
of edges with infectious (resp., removed, susceptible) alter among those
having susceptible ego. We also introduce
%
%
\begin{equation}\label{deftheta}
\theta_t = \exp\biggl(-r \int_0^t \bpI_s \,\rmd {s}\biggr),
\end{equation}
the probability that a degree one node remains susceptible until time
$t$. For
any $f \in\mathcal B_b(\N)$,
%
%
\begin{eqnarray}
\label{limitereseauinfiniS}
\langle\bmuS_t, f\rangle&=& \sum_{k\in\N} \bmuS_0(k) \theta^k_t
f(k),\\
\label{limitereseauinfiniSI}
\langle\bmuSI_t, f\rangle&=& \langle
\bmuSI_0, f\rangle-\int_0^t \beta\langle
\bmuSI_s, f\rangle\,\rmd {s} \nonumber\\
&&{}
+ \int_0^t \sum_{k\in\N} rk \bpI_s \mathop{\sum_{j, \ell,
m\in
\N}}_{j+\ell+m=k-1}
\pmatrix{k-1\cr j,\ell,m}
(\bpI_s)^{j}(\bpR_s)^{\ell}(\bpS_s)^{m}\nonumber\\
&&\qquad\hspace*{143pt}{}\times f(m)\bmuS_s(k)\,\rmd {s}\\
&&{} + \int_0^t \sum_{k\in\N} rk \bpI_s \bigl(1+(k-1) \bpI_s\bigr) \nonumber\\
&&\qquad\quad\hspace*{8.2pt}{}\times\sum
_{k'\in
\N^*} \bigl(f(k'-1)-f(k')\bigr) \frac{k'\bmuSI_s(k')}{\bNIS_s} \bmuS
_s(k)\,\rmd {s},
\nonumber
\\
\label{limitereseauinfiniSR}
\langle\bmuSR_t, f\rangle&=& \langle \bmuSR_0, f\rangle + \int_0^t \beta\langle\bmuSI_s,
f\rangle\,\rmd {s} \nonumber\\
&&{} + \int_0^t \sum_{k\in\N} rk\bpI_s (k-1) \bpR_s \sum_{k'\in
\N^*} \bigl(f(k'-1)-f(k')\bigr) \\
&&\qquad\quad\hspace*{98.7pt}{}\times\frac{k'\bmuRS_s(k')}{\bNRS_s} \bmuS
_s(k)\,\rmd {s}.\nonumber
\end{eqnarray}

We denote by $\bar{S}_t$ (resp., $\bar{I}_t$ and $\bar{R}_t$) the
mass of the measure $\bmuS_t$ (resp.,~$\bmuSI_t$ and $\bmuSR_t$). As
for the finite graph, $\bmuS_t/\bar{S}_t$
(resp., $\bmuSI_t/\bar{I}_t$ and $\bmuSR_t/\bar{R}_t$) is the
probability degree distribution
of the susceptible individuals (resp., the probability distribution of
the degrees of the infectious and removed individuals toward the
susceptible ones).

Let us give a heuristic explanation of equations (\ref
{limitereseauinfiniS})--(\ref{limitereseauinfiniSR}). Remark that the
graph in the limit is infinite. The probability that an individual of
degree $k$ has been infected by none of her $k$ edges is $\theta_t^k$
and equation (\ref{limitereseauinfiniS}) follows. In (\ref
{limitereseauinfiniSI}), the first integral corresponds to infectious
individuals being removed. In the second integral, $rk\bpI_s$ is the
rate of infection of a given susceptible individual of degree~$k$. Once
she gets infected, the multinomial term determines the number of edges
connected to susceptible, infectious and removed neighbors. Multi-edges
do not occur. Each infectious neighbor has a degree chosen in the
size-biased distribution $k' \bmuSI(k')/\bNIS$ and the number of
edges to $\cS_t$ is reduced by 1. This explains the third integral.
Similar arguments hold for (\ref{limitereseauinfiniSR}).

Choosing $f(k)=\ind_{i}(k)$, we obtain the
following countable system of ordinary differential equations (ODEs):
%
%
\begin{eqnarray}\label{systemelong}\hspace*{18pt}
\bar{\mu}_t^{\sscS}(i)&=& \bar{\mu}^{\sscS}_0(i) \theta
^i_t,\nonumber\\
\bmuSI_t(i) &=& \bmuSI_0(i)
+ \int_0^t \biggl\{ r \bpI_s \sum_{j,\ell\geq0} (i+j+\ell+1)\bar
{\mu
}_s^{\sscS}(i+j+\ell+1) \nonumber\\
&&\qquad\quad\hspace*{62.3pt}{}\times\pmatrix{i+j+\ell\cr i,j,\ell} (\bpS
_s)^i(\bpI
_s)^j (\bpR_s)^{\ell}\nonumber\\
&&\hspace*{56.1pt}{} + \bigl(r(\bpI_s)^2\langle\bar{\mu}_s^{\sscS},\chi^2-\chi\rangle+
r\bpI_s \langle\bar{\mu}^{\sscS}_s,\chi\rangle\bigr)\\
&&\qquad\quad\hspace*{34pt}{}\times
\frac{(i+1)\bmuSI_s(i+1)-i\bmuSI_s(i)}{\langle\bmuSI_s,\chi
\rangle}-\beta\bmuSI_s(i)\biggr\} \,\rmd {s},\nonumber\\
\bmuSR_t(i) &=& \bmuSR_0(i)+\int_0^t \biggl\{ \beta
\bmuSI_s(i)\nonumber\\
&&\hphantom{\int_0^t \biggl\{\bmuSR_0(i)\,+ }{}+ r \bpI_s \langle
\bar{\mu}_s^{\sscS},\chi^2-\chi\rangle\bar{p}^{\sRr}_s
\frac{(i+1)\bmuSR_s(i+1)-i\bmuSR_s(i)}{\langle
\bmuSR_s,\chi\rangle}\biggr\} \,\rmd {s}.\nonumber
\end{eqnarray}
It is noteworthy to say that this system is similar but not identical
to that in Ball and Neal \cite{ballneal}.
Our equations differ since our mechanism is not the same (compare
Section \ref{sectiondynamics} with Section 5 in \cite{ballneal}). We
emphasize that the number of links of an individual to $\scS$ decreases
as the epidemic progresses, which modifies her infectivity.

The system (\ref{limitereseauinfiniS})--(\ref{limitereseauinfiniSR})
allows us to recover the equations proposed by Volz \cite{volz}, Table
3, page 297. More
precisely, the dynamics of the epidemic is obtained by solving the
following closed system of four ODEs, referred to as
Volz's equations in the sequel. The latter are obtained directly from
(\ref{limitereseauinfiniS})--(\ref{limitereseauinfiniSR}) and the
definitions of $\bar{S}_t$, $\bar{I}_t$, $\bpI_t$ and $\bpS_t$
which relate these quantities to the measures $\bmuS_t$ and $\bmuSI
_t$. Let
%
%
\begin{equation}
g(z)=\sum_{k\in
\N}\bar{\mu}_0^{\sscS}(k) z^k
\end{equation}
be the generating function for the
initial degree distribution of the susceptible individuals
$\bar{\mu}_0^{\sscS}$, and let $\theta_t= \exp(-r\int_0^t \bar
{p}_s^{{\II}} \,\rmd {s})$. Then, the epidemic can be approximated by
the solution of the four following
ODEs:
%
%
\begin{eqnarray}
\label{volz1}
\bar{S}_t &=& \langle\bar{\mu}_t^{\sscS},\mathbf1\rangle=
g(\theta
_t),\\
\label{volz2}
\bar{I}_t &=& \langle\bmuIS_t,\mathbf1\rangle= \bar{I}_0+\int_0^t
\bigl( r\bpI_s \theta_s g'(\theta_s)-\beta\bar{I}_s\bigr) \,\rmd {s},\\
\label{volz3}
\bpI_t &=& \bpI_0
+ \int_0^t \biggl( r \bpI_s \bpS_s\theta_s\frac{ g''(\theta
_s)}{g'(\theta_s)}-r \bpI_s(1-\bpI_s) -\beta\bpI_s\biggr)\,\rmd {s},\\
\label{volz4}
\bpS_t &=& \bpS_0+\int_0^t r \bpI_s \bpS_s
\biggl(1-\theta_s\frac{g''(\theta_s)}{g'(\theta_s)}\biggr)\,\rmd {s}.
\end{eqnarray}
Here, the\vspace*{1pt} graph structure appears through the generating
function $g$. In~(\ref{volz2}), we
see that the classical contamination terms $r \bar{S}_t \bar{I}_t$
(mass action) or $r\bar{S}_t\bar{I}_t/\allowbreak(\bar{S}_t+\bar{I}_t)$
(frequency dependence) of mixing SIR models (e.g., \cite
{andersonbritton,arazozaclemencontran}) are replaced by $r\bpI_t
\theta_t
g'(\theta_t)=r\bNIS_t$. The fact that new infectious individuals are
chosen in the size-biased distribution is hidden in the term
$g''(\theta_t)/g'(\theta_t)$.

The beginning of the epidemic and computation of the reproduction
number, when the numbers of infected individuals and of contaminating
edges are small and when Volz's deterministic approximation does not
hold, makes the object of another study.

\section{SIR model on a configuration model graph}\label{sectionmodel}

In this section, we introduce configuration model graphs and describe
the propagation of SIR on such graphs.

\subsection{Configuration model graph}
\label{secconfigurations-graph}

Graphs at large can be mathematically represented as matrices
with integer entries: to each graph corresponds an adjacency
matrix, the $(x, y)$th coefficient of which is the number of edges
between the vertices $x$ and $y$. Defining the distribution of a random
graph thus amounts to choosing a sigma-field and a probability measure
on the space~$\N^{\N^*\times\N^*}$, where $\N^*=\N\setminus\{0\}
$. Another approach is to construct
a random graph by modifying progressively a given graph, as in Erd\"
{o}s--Renyi model. Several other constructions are possible such as the
preferential
attachment model, the threshold graphs, etc. 

Here, we are interested in the configuration model (CM) proposed by
Bollob\'{a}s~\cite{bollobas2001}, Molloy and Reed \cite{molloyreed}
(see also \cite{newmanSIAM,newmanstrogatzwatts,durrett,vanderhofstad})
and which models graphs with specified degree distribution and
independence between the degrees of neighbors. As shown by statistical
tests, these models might be realistic in describing community
networks. See, for instance, Cl\'{e}men\c{c}on et al. \cite
{clemencondearazozarossitran} for dealing with the spread of the
HIV--AIDS disease among the homosexual community in Cuba.

We recall its construction (see, e.g., \cite
{durrett,vanderhofstad}). Suppose we are given the number of vertices,
$N$ and i.i.d. random variables (r.v.)
$d_1,\ldots, d_N$ with distribution $(p_k)_{k\in\N}$ that represent
the degrees of each vertex. To the vertex~$i$ are associated $d_i$
half-edges. To construct an edge, one chooses two open half-edges
uniformly at random and pair them together.

Remark that this linkage procedure does not exclude self-loops or
multiple edges. In the following, we are interested in a large number
of nodes with a fixed degree
distribution, hence self-loops and multiple edges
become less and less apparent in the global picture (see, e.g.,
\cite{durrett}, Theorem 3.1.2).

Notice that the condition for the existence with positive probability
of a~giant component
is that the expectation of the size biased distribution is larger than 1:
\[
\sum_{k\in\N}(k-1)\frac{kp_k}{\sum_{\ell\in\N} \ell
p_\ell}>1.
\]
This is connected with the fact that the Galton--Watson tree with this
offspring distribution is supercritical (see \cite{durrett}, Section
3.2, page 75, for details).

\subsection{SIR epidemic on a CM graph}\label{sectiondynamics}

\begin{figure}

\includegraphics{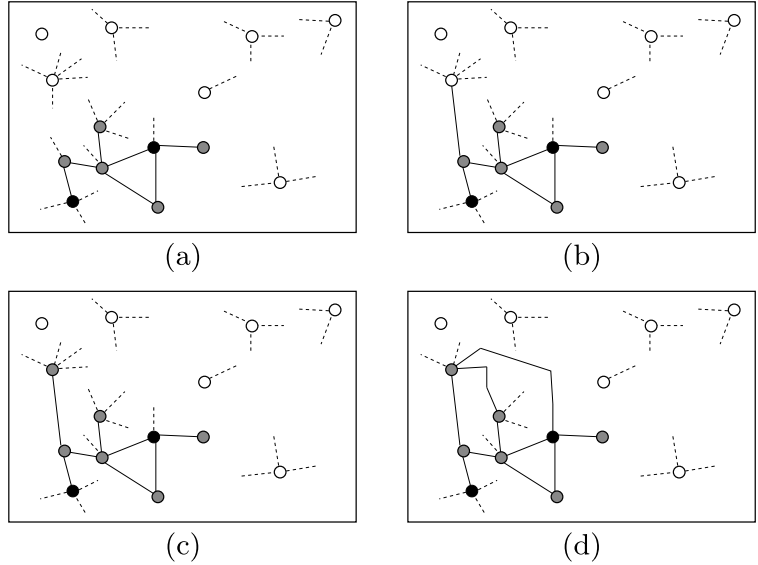}

\caption{Infection process. Arrows provide the infection tree.
Susceptible, infectious and removed individuals are colored in white,
grey and black, respectively. \textup{(a)} The degree of each individual is
known, and for each infectious (resp., removed) individual, we know her
number of edges of type $\cI\cS$ (resp., $\cRr\cS$). \textup{(b)}, \textup{(c)} A
contaminating half-edge is chosen, and say that a susceptible of degree
$k$ is infected at time $t$. 
The contaminating edge is drawn in bold line. 
\textup{(d)} Once the susceptible individual has been infected, we determine how
many of her remaining arrows are linked to the classes $\cI$ and $\cRr
$. If we denote by $j$ and $\ell$ these numbers, then $\NIS_t=\NIS
_{t_-}-1+(k-1)-j-\ell$ and $\NRS_t=\NRS_{t_-}-\ell$.}
\label{figdyn2}
\end{figure}

We now propagate an epidemic on a CM graph of size $N$
(see Figure \ref{figdyn2}).
The disease can be transmitted from infectious nodes to neighboring
susceptible nodes and removed nodes cannot be reinfected.

%

Suppose that at initial time, we are given a set of \textit{susceptible}
and \textit{infectious} nodes together with their degrees. The graph of
relationships between the individuals is, in fact, irrelevant for
studying the propagation of the disease. The minimal information
consists in the sizes of the classes $\scS$, $\I$, $\Rr$ and the
number of edges to the class $\scS$ for every infectious or removed
node. Thus, each node of class $\cS$ comes with a
given number of half-edges of undetermined types; each node of class
$\I$ (resp., $\Rr$) comes with
a number of $\cI\cS$ (resp., $\Rr\scS$) edges. The numbers of $\I
\Rr
$, $\I\I$ and $\Rr\Rr$ edges need not to be retained.

The evolution of the SIR epidemic on a CM graph can be described as
follows. To each $\cI\cS$-type half-edge is associated an independent
exponential clock with parameter $r$ and
to each $\cI$ vertex is associated an independent exponential clock
with parameter $\beta$. The first of all these clocks that rings
determines the next event.
\begin{enumerate}[\textit{Case} 2.]
\item[\textit{Case} 1.] If the clock that rings is associated to an $\cI$
individual, the latter recovers. Change her status from $\cI$ to $\cRr
$ and
the type of her emanating half-edges accordingly: $\cI\cS$ half-edges
become $\cRr\cS$ half-edges.
\item[\textit{Case} 2.] If the clock is associated with a half $\cI\cS$-edge,
an infection occurs.
\begin{enumerate}[\textit{Step} 2.]
\item[\textit{Step} 1.] Match randomly the $\cI\cS$-half-edge that has rung
to a half-edge belonging to a susceptible.
\item[\textit{Step} 2.] This susceptible is the newly infected. Let $k$ be her degree.
Choose\vadjust{\goodbreak} uniformly $k-1$ half-edges among all the available half-edges
(they either are of type $\cI\cS$, $\cRr\cS$ or emanate from~$\scS$). Let $m$ (resp., $j$ and $\ell$) be the number of $\scS\scS$-type
(resp., of $\cI\cS$ and of $\cRr\cS$-type) half-edges drawn among
these $k-1$ half-edges.
\item[\textit{Step} 3.] 
The chosen half-edges of type $\cI\cS$ and $\cRr\cS$ determine the
infectious or removed neighbors of the newly infected individual. The
remaining $m$ edges of type $\cS\cS$ remain open in the sense that the
susceptible neighbor is not fixed. Change the status of the $m$ (resp.,
$j$, $\ell$) $\cS\cS$-type (resp., $\cI\cS$-type, $\cRr\cS$-type) edges
created to $\cS\cI$-type (resp., $\cI\cI$-type, $\cRr\cI$-type).
\item[\textit{Step} 4.] Change the status of the newly infected from $\cS$ to
$\cI$.
\end{enumerate}
\end{enumerate}
We then wait for another clock to ring and repeat the procedure.

We only need three descriptors of the system to obtain a Markovian
evolution, namely, the three degree distributions introduced in (\ref
{muS}).

For a measure $\mu\in\mathcal{M}_F(\N)$, we denote by $F_\mu
(m)=\mu(\{0,\ldots,m\})$, $m\in\N$, its cumulative distribution
function. We introduce $F^{-1}_\mu$ its right inverse (see the
\hyperref[appendiceA]{Appendix}). Then, for all $0\leq i\leq S_t$
(resp., $0\leq i \leq I_t$ and $0\leq i\leq R_t$),
\[
\gamma_i(\muS_t)=F_{\muS_t}^{-1}(i)\qquad [\mbox{resp., }\gamma_i(\muSI
_t)=F_{\muSI_t}^{-1}(i), \gamma_i(\muSR_t)=F_{\muSR_t}^{-1}(i)]
\]
represents the degree at $t$ of the $i$th susceptible individual
(resp., the number of edges to $\cS$ of the $i$th infectious
individual and of the $i$th removed individual) when individuals are
ranked by increasing degrees (resp., by number of edges to $\cS$).
\begin{ex}\label{exempleFmu}
Consider, for instance, the measure $\mu=2\delta_1 + 3\delta
_5+\delta_7$.
Then, the atoms $1$ and $2$ are at level $1$, the atoms $3$, $4$ and
$5$ are at level~$5$ and the atom $6$ is at level $7$. We then have that
$\gamma_1(\mu)=F_{\mu}^{-1}(1)=1$, $\gamma_2(\mu)=1$, $\gamma
_3(\mu)=\gamma_4(\mu)=\gamma_5(\mu)=5$ and $\gamma_6(\mu)=7$.
\end{ex}

From $t$, and because of the properties of
exponential distributions, the next event will take place in a time
exponentially
distributed with parameter $r\NIS_t+\beta I_t$. Let $T$
denote the time of this event.
\begin{longlist}[\textit{Case} 1.]
\item[\textit{Case} 1.] The event corresponds to a removal, that is, a
node goes from status $\cI$ to status $\cRr$.
Choose uniformly an integer $i$ in $I_{T^-}$, then
update the measures $\muSI_{T_-}$ and $\muSR_{T_-}$:
\[
\muSI_T=\muSI_{T^-} - \delta_{\gamma_i(\muSI_{T_-})}
\quad\mbox{and}\quad
\muSR_T=\muSR_{T^-}+\delta_{\gamma_i(\muSI_{T_-})}.
\]
The probability that a given integer $i$ is drawn is
$1/I_{T^-}$. 

\item[\textit{Case} 2.] The event corresponds to a new infection.
We choose uniformly a~half-edge with susceptible ego, and this ego
becomes infectious. The global rate of infection is $r\NIS_{T_-}$
and the probability of choosing a susceptible individual of degree $k$
for the new infectious is $k\muS_{T_-}(k)/\NS_{T_-}$.
We define by
%
%
\begin{equation}
\label{eqdeflambda}
\lambda_{T_-}(k)=rk\frac{\NIS_s}{\NS_s}
\end{equation}
the rate of infection of a \textit{given} susceptible of degree $k$ at
time $T_-$. This notation was also used in Volz \cite{volz}.

%
The newly infective may have several links with infectious or removed
individuals. The probability, given that the degree of the individual
is $k$ and that
$j$ (resp.,~$\ell$) out of her $k-1$ other half-edges (all but the
contaminating $\cI\cS$ edge) are chosen to be
of type $\cI\cI$ (resp., $\cI\cRr$), according to Step 2, is given
by the following multivariate hypergeometric distribution:
%
%
\begin{equation}
\label{eqdefp}
p_{T_-}(j,\ell| k-1)=\frac{{\NIS_{T_-}-1 \choose j} {\NRS_{T_-}
\choose\ell} {\NS_{T_-}-\NIS_{T_-}-\NRS_{T_-} \choose k-1-j-\ell
}}{{\NS_{T_-}-1 \choose k-1}}\cdotp
\end{equation}
%
Finally, to update the values of $\muSI_T$ and $\muSR_T$ given $k$,
$j$ and $\ell$, we have to choose the infectious and removed
individuals to which the newly infectious is linked;\vadjust{\goodbreak} some of their
edges, which were $\cI\cS$ or $\cRr\cS$, now become $\cI\cI$ or
$\cRr\cI$.
We draw two sequences $u=(u_1,\ldots,u_{I_{T_-}})$ and $v=(v_1,\ldots
,v_{R_{T_-}})$ that will indicate how many links each infectious or
removed individual has to the newly contaminated individual. There
exist constraints on $u$ and $v$: the number of edges recorded by $u$
and $v$ cannot exceed the number of existing edges.
Let us define the set
%
%
\begin{equation}
\label{eqdefU}
\mathcal{U}=\bigcup_{m=1}^{+\infty} \N^m,
\end{equation}
and for all finite integer-valued measure $\mu$ on $\N$, and all
integer $n \in\N$, we define the subset
%
%
\begin{eqnarray}
\label{eqdefU2}
&&\mathcal U(n,\mu)=\Biggl\{u=\bigl(u_1,\ldots,u_{\langle\mu,\car\rangle}\bigr) \in
\mathcal U \mbox{ such that }\nonumber\\[-8pt]\\[-8pt]
&&\qquad\hspace*{33.7pt}\forall i\in\{1,\ldots,\langle\mu,\car\rangle\}, u_i\leq
F^{-1}_\mu(i)\mbox{ and }\sum_{i=1}^{\langle\mu,\mathbf1 \rangle
} u_i=n\Biggr\}.
\nonumber
\end{eqnarray}
Each sequence $u\in\mathcal{U}(n,\mu)$ provides a possible
configuration of how the $n$ connections of a given individual can be
shared between neighbors whose degrees are summed up by $\mu$. The
component $u_i$, for $1\leq i\leq\langle\mu, 1\rangle$, provides
the number of edges that this individual shares with the $i$th
individual. This number is necessarily smaller than the degree $\gamma
_i(\mu)=F_{\mu}^{-1}(i)$ of individual~$i$. Moreover, the $u_i$'s sum
to~$n$. The probabilities of the draws of $u$ and $v$ that provide,
respectively, the number of edges
$\cI\cS$ which become $\cI\cI$ per infectious individual and the
number of edges $\cRr\cS$ which become $\cRr\cI$ per removed
individual are given by
%
%
\begin{eqnarray}\label{etape14}
\rho(u | j+1,\muSI_{T_-}) &=& \frac{\prod_{i=1}^{I_{T_-}}{\gamma
_i(\muSI_{T_-}) \choose u_i}}{{ \NIS_{T_-} \choose j+1}}\ind_{u\in
\mathcal{U}(j+1,\muSI_{T_-})},\nonumber\\[-8pt]\\[-8pt]
\rho(v | \ell,\muSR_{T_-}) &=& \frac{\prod_{i=1}^{R_{T_-}}{ \gamma
_i(\muSR_{T_-}) \choose v_i}}{{ \NRS_{T_-} \choose\ell}}\ind_{v\in
\mathcal{U}(\ell,\muSR_{T_-})}.\nonumber
\end{eqnarray}
Then, we update the measures as follows:
%
%
\begin{eqnarray}\label{etape1546}
\muS_T&=&\mu^{\sscS}_{T^-}-\delta_k, \nonumber\\
\muSI_T&=&\muSI_{T^-} +\delta_{k-1-j-\ell}+\sum_{i=1}^{{\II}_{T_-}}
\delta_{\gamma_{i}(\muSI_{T_-})-u_{i}}-\delta_{\gamma_{i}(\muSI
_{T_-})},\\
\muSR_T&=&\muSR_{T^-} +\sum_{i'=1}^{\sRr_{T_-}} \delta_{\gamma
_{i'}(\muSR_{T_-})-v_{i'}}-\delta_{\gamma_{i'}(\muSR
_{T_-})}.\nonumber
\end{eqnarray}
\end{longlist}





\subsection{Stochastic differential equations}


Here, we propose stochastic differential equations (SDEs) driven
by\vadjust{\goodbreak}
Poisson point measures (PPMs) to describe the evolution of the degree
distributions (\ref{muS}), following the inspiration of \cite
{arazozaclemencontran,fourniermeleard}. We consider two PPMs: $\rmd
{Q}^1(s,k,\theta_1,j,\ell,\theta_2, u,\theta_3,v,\theta_4)$ and
$\rmd {Q}^2(s,i)$ on $\R_+\times E_1$ with $E_1:= \N\times\R
_+ \times\N\times\N\times\R_+ \times\U\times\R_+\times\U
\times\R_+ $ and $\R_+\times\N$ with intensity measures $\rmd
{q}^1(s,k,\theta_1,j,\ell,\theta_2,u,\theta_3,v,\theta
_4)=\rmd {s} \otimes\rmd {n}(k)\otimes\rmd \theta
_1 \otimes\rmd {n}(j)\otimes\rmd {n}(\ell)\otimes
\rmd \theta_2\otimes\rmd {n}(u)\otimes\rmd \theta
_3 \otimes\rmd {n}(v)\otimes\rmd \theta_4$ and $
\rmd{q}^2(s,i)=\beta\,\rmd {s}\otimes\rmd {n}(i)$, where
$\rmd {s}$, $\rmd \theta_1$, $\rmd \theta_2$,
$\rmd \theta_3$ and $\rmd \theta_4$ are Lebesgue measures
on $\R_+$, where $\rmd {n}(k)$, $\rmd {n}(j)$, $
\rmd{n}(\ell)$ are counting measures on $\N$ and where $
\rmd{n}(u)$, $\rmd {n}(v)$ are counting measures on $\U$.

The point measure $Q^1$ provides possible times at which an infection
may occur. Each of its atoms is associated with a possible infection
time $s$, an integer $k$ which gives the degree of the susceptible
being possibly infected, the number $j+1$ and $\ell$ of edges that
this individual has to the sets $\cI_{s_-}$ and~$\cRr_{s_-}$. The
marks $u$ and $v\in\U$ are as in the previous section. The
marks~$\theta_1$, $\theta_2$ and~$\theta_3$ are auxiliary variables used
for the construction [see (\ref{eqmuI}) and~(\ref{eqmuR})].

The point measure $Q^2$ gives possible removal times. To each of its
atoms is associated a possible removal time $s$ and the number $i$ of
the individual that may be removed.

The following SDEs describe the evolution of the epidemic: for all
$t\ge0$,
%
%
\begin{eqnarray}
\label{eqmuS}\qquad\quad
\muS_t &=& \muS_0 - \int_0^t\int_{E_1} \delta_{k}\ind_{\theta
_1\leq\lambda_{s_-}(k)\muS_{s_-}(k)}\nonumber\\[-9.5pt]\\[-9.5pt]
&&\hphantom{\muS_0 - \int_0^t\int_{E_1}}{}\times\ind_{\theta_2\leq p_{s_-}(j,\ell|k-1)} \ind
_{\theta_3\leq\rho(u|j+1,\muSI_{s_-})} \ind_{\theta_4\leq\rho
(v|\ell,\muSR_{s_-})}
\,\rmd {Q}^1, \nonumber\\[-2pt]
\label{eqmuI}
\muSI_t &=& \muSI_0+\int_0^t\int_{E_1} \Biggl(\delta_{k-(j+1+\ell
)}+\sum
_{i=1}^{I_{s_-}}\bigl(\delta_{\gamma_i(\muSI_{s_-})-u_i}-\delta_{\gamma
_i(\muSI_{s_-})}\bigr)\Biggr)\nonumber\\[-2pt]
&&\hphantom{\muSI_0+\int_0^t\int_{E_1}}{} \times\ind_{\theta_1\leq\lambda_{s_-}(k)\mu
_{s_-}^{\sscS}(k)}\ind_{\theta_2\leq p_{s_-}(j,\ell|k-1)}
\ind
_{\theta_3\leq\rho(u|j+1,\muSI_{s_-})}\nonumber\\[-9.5pt]\\[-9.5pt]
&&\hphantom{\muSI_0+\int_0^t\int_{E_1}}{} \times
\ind_{\theta_4\leq\rho(v|\ell,\muSR_{s_-})}
\,\rmd {Q}^1 \nonumber\\[-2pt]
&&{} - \int_0^t\int_{\N} \delta_{\gamma_i(\muSI
_{s_-})} \ind_{i\leq I_{s_-}}\,\rmd {Q}^2, \nonumber\\[-2pt]
\label{eqmuR}
\muSR_t &=& \muSR_0+\int_0^t\int_{E_1} \Biggl(\sum_{i=1}^{R_{s_-}} \bigl(\delta
_{\gamma_{i}(\muSR_{s_-})-v_{i}}-\delta_{\gamma_{i}(\muSR
_{s_-})}\bigr)\Biggr)\nonumber\\[-2pt]
&&\hphantom{\muSR_0+\int_0^t\int_{E_1} }{} \times\ind_{\theta_1\leq\lambda_{s_-}(k) \mu
^{\sscS
}_{s_-}(k)}\ind_{\theta_2\leq
p_{s_-}(j,\ell|k-1)}
\ind_{\theta_3\leq\rho(u|j+1,\muSI_{s_-})}\nonumber\\[-9.5pt]\\[-9.5pt]
&&\hphantom{\muSR_0+\int_0^t\int_{E_1} }{} \times
\ind_{\theta_4\leq\rho(v|\ell,\muSR_{s_-})} \,\rmd {Q}^1\nonumber\\[-2pt]
&&{} + \int_0^t\int_{\N} \delta_{\gamma_i(\muSI
_{s_-})} \ind_{i\leq I_{s_-}}\,\rmd {Q}^2,\nonumber
\end{eqnarray}
where we write $\rmd {Q}^1$ and $\rmd {Q}^2$ instead of
$\rmd {Q}^1(s,k,\theta_1,j,\ell,\theta_2, u,\theta_3,v,\theta
_4)$ and $\rmd {Q}^2(s,i)$
to simplify the notation.
\begin{prop}For any given initial conditions $\mu^{\sscS}_0$, $\mu
^{\sscS
{\II}}_0$ and $\muSR_0$ that are integer-valued measures on $\N$ and for
PPMs $Q^1$ and $Q^2$, there exists a~unique strong
solution to the SDEs (\ref{eqmuS})--(\ref{eqmuR}) in the space
$\D(\R_+,\allowbreak(\mathcal{M}_F(\N))^3)$, the Skohorod space of c\`{a}dl\`
{a}g functions with values in $(\mathcal{M}_F(\N))^3$.
\end{prop}
\begin{pf}For the proof we notice that for every $t\in\R_+$, the
measure $\mu^{\sscS}_t$ is dominated by $\mu^{\sscS}_0$ and the measures
$\mu^{{\II}\sscS}_t$ and $\mu^{\sRr\sscS}_t$ have\vspace*{1pt} a mass bounded by
$\langle\mu^{\sscS}_0+\mu_0^{{\II}\sscS}+\muSR_0,1\rangle$ and a support
included in
$[ [
0,\max\{\max(\supp(\mu^{\sscS}_0))$, $\max(\supp(\mu^{{\II}\sscS
}_0))$, $\max(\supp(\muSR_0))\}]
]$. The
result then follows the steps of~\cite{fourniermeleard} and
\cite{chithese} (Proposition~2.2.6).
\end{pf}

The course of the epidemic can be deduced from (\ref{eqmuS}), (\ref
{eqmuI}) and (\ref{eqmuR}). For the sizes $(S_t,I_t,R_t)_{t\in\R_+}$
of the different classes, for instance, we have with the choice of
$f\equiv1$, that for all $t\ge0$, $S_t=\langle\mu^{\sscS
}_t,\mathbf
{1} \rangle$, $ I_t= \langle\muSI_t,\mathbf1\rangle$ and $ R_t=
\langle\muSR_t,\mathbf1\rangle$. Writing the semimartingale
decomposition that results from standard stochastic calculus for jump
processes and SDE driven by PPMs (e.g.,
\cite{fourniermeleard,ikedawatanabe,jacod}), we obtain, for example,
%
%
\begin{equation}\label{eqIt}
I_t= \langle\muSI_t,\mathbf1\rangle= I_0+\int_0^t\biggl( \sum_{k\in\N}
\mu_s^{\sscS}(k)\lambda_s(k)
- \beta I_s\biggr) \,\rmd {s}+M^{{\II}}_t,
\end{equation}
where $M^{{\II}}$ is a square-integrable martingale that can be written
explicitly with the compensated PPMs of $Q^1$ and $Q^2$, and with
predictable quadratic variation given for all $t\ge0$ by
\[
\langle M^{{\II}}\rangle_t=\int_0^t \sum_{k\in\N} \bigl(\mu_s^{\sscS
}(k)\lambda_s(k)+\beta I_s\bigr)\,\rmd {s}.
\]
Another quantities of interest are the numbers of edges of
the different types~$N^{\sscS}_t$, $\NIS_t$, $\NRS_t$. The latter
appear as the first moments of the measures $\muS_t$, $\muSI_t$ and
$\muSR_t\dvtx N^{\sscS}_t= \langle\muS_t,\chi\rangle$, $\NIS_t=
\langle\muSI_t,\chi\rangle$ and $\NRS_t= \langle\muSR_t,\chi
\rangle$.

\section{Large graph limit}\label{sectionlargegraph}

Volz \cite{volz} proposed a parsimonious
deterministic approximation to describe the epidemic dynamics when
the population is large. However, the stochastic processes are
not clearly defined and the convergence of the SDEs to the four ODEs
that Volz proposes is stated but not proved. Using the construction
that we developed in Section \ref{sectiondynamics}, we provide
mathematical proofs of Volz's equation, starting from a finite graph
and taking the limit when the size of the graph tends to infinity.
Moreover, we see that the three distributions
$\mu^{\sscS}$, $\muSI$ and $\muSR$ are at the core of the problem and
encapsulate the evolution of the process.

\subsection{Law of large numbers scaling}\label{sectionrenorm}

We consider sequences of measures $(\mu^{n,\sscS})_{n\in\N}$,
$(\mu^{n,{\II}\sscS})_{n\in\N}$ and $(\mu^{n,\sRr\sscS})_{n\in
\N}$ such
that for any $n\in{\mathbb N}^*$, $\mu^{n,\sscS}$, $\mu^{n,{\II
}\sscS}$ and
$\mu^{n,\sRr\sscS}$ satisfy (\ref{eqmuS})--(\ref{eqmuR}) with initial
conditions $\mu^{n,\sscS}_0$, $\mu^{n,{\II}\sscS}_0$ and $\mu
^{n,\sRr\sscS}_0$.
We denote by
$\scS_t^n$, $\I_t^n$ and $\Rr^n_t$ the subclasses of
susceptible,\vspace*{1pt}
infectious or removed individuals at time $t$, and by $\NSnn_t$,
$\NISnn_t$ and $\NRSnn_t$, the number of edges with susceptible ego,
infectious ego and susceptible alter, removed ego and susceptible alter.
The number of vertices of each class are denoted
$I^n_t$, $S^n_t$ and~$R^n_t$.
The total size of the population is finite and equal to
$S^n_0+I^n_0+R^n_0$. The size of the population and the number of edges tend
to infinity proportionally to $n$.

We scale the measures the following way. For any $n \ge0$, we set
\[
\muSIn_t=\frac{1}{n}\mu^{n,{\II}\sscS}_t
\]
for all $t\geq0$ (and accordingly, $\muSn_t$ and $\muSRn_t$).
Then, we denote
\[
\NISn_t=\bigl\langle\muISn,\chi \bigr\rangle=\frac{1}{n} \NISnn_t
\quad\mbox{and}\quad
\In_t=\bigl\langle\muISn_t,\mathbf1 \bigr\rangle=\frac{1}{n} I^n_t
\]
and accordingly, $\NSn_t$, $\NRSn_t$, $\Sn_t$ and $\Rn_t$.

We assume that the initial conditions satisfy:
\begin{hyp}\label{hypconvcondinit}
$\!\!\!$The sequences $(\muSn_0)_{n\in \N }$, $(\muSIn_0)_{n\in \N }$ and
$(\muSRn_0)_{n\in \N}$ converge to measures $\bmuS_0$, $\bmuSI_0$ and
$\bmuSR_0$ in $\mathcal{M}_F(\N)$ embedded with the weak convergence
topology.
\end{hyp}
\begin{Rque}\label{rqueconvcondinit}
(1) Assumption \ref{hypconvcondinit} entails that the initial
(susceptible and infectious) population size is of order $n$
if $\bmuS_0$ and $\bmuSI_0$ are nontrivial.


(2) In case the distributions underlying the measures
$\mu^{n,\sscS}_0$, $\mu^{n,{\II}\sscS}_0$ and $\mu^{n,\cRSr}_0$
do not depend on the total number of vertices (e.g., Poisson,
power-laws or geometric distributions), Assumption
\ref{hypconvcondinit} can be viewed as a law of large numbers. When the
distributions depend on the total number of vertices~$N$ (as in
Erd\"{o}s--Renyi graphs), there may be scalings under which
Assumption~\ref{hypconvcondinit} holds. For Erd\"{o}s--Renyi graphs, for instance,
if the probability $\rho_N$ of connecting two vertices satisfies
$\lim_{N\rightarrow +\infty} N\rho_N=\lambda$, then we obtain in the
limit a Poisson distribution with parameter $\lambda$.

(3) In (\ref{limitereseauinfiniSI}), notice the
appearance of the size biased degree distribution
$k\bmuS_s(k)/N^{\sscS}_s$. The latter reflects the fact
that, in the CM, individuals having large degrees have higher
probability to connect than individuals having small
degrees. Thus, there is no reason why the degree distributions of
the susceptible individuals $\bmuS_0/\bar{S}_0$
and the distribution $\sum_{k\in\N}p_k\delta_k$ underlying the CM
should coincide. Assumption \ref{hypconvcondinit} tells\vspace*{1pt} us indeed
that the initial infectious population size is of order~$n$. Even
if $\bar{I}_0 /\bar{S}_0$ is very small, the biased distributions
that appear imply that the degree distribution
$\bmuSI_0/\bar{I}_0$ should
have a larger expectation than the degree distribution
$\bmuS_0/\bar{S}_0$.
\end{Rque}

We obtain rescaled SDEs which are the same as the SDEs (\ref
{eqmuS})--(\ref{eqmuR}) parameterized by $n$. For all $t\ge0$,
%
%
\begin{eqnarray}
\label{eqmuSrenorm}
\muSn_t &=& \muSn_0-\frac{1}{n}\int_0^t\int_{E_1} \delta_{k}\ind
_{\theta_1\leq\lambda^n_{s_-}(k)n\muSn_{s_-}(k)}
\ind_{\theta_2\leq p^n_{s_-}(j,\ell|k-1)} \nonumber\\[-8pt]\\[-8pt]
&&\hspace*{77pt}{} \times\ind_{\theta_3\leq\rho(u|j+1,n\muISn
_{s_-})} \ind_{\theta_4\leq\rho(v|\ell,n\muRSn_{s_-})}
\,\rmd {Q}^1,\nonumber\\
\label{eqmuIrenorm}
\muISn_t &=& \muSIn_0+\frac{1}{n}\int_0^t\!\int_{E_1}\!\Biggl(\!\delta
_{k-(j+1+\ell)}+\sum_{i=1}^{I^n_{s_-}}\bigl(\delta_{\gamma_{i}(n\muISn
_{s_-})-u_i}-\delta_{\gamma_{i}(n\muISn_{s_-})}\bigr)\!\Biggr)\hspace*{-37pt}\nonumber\\
&&\hspace*{79.3pt}{} \times\ind_{\theta_1\leq\lambda^n_{s_-}(k)n\muSn
_{s_-}(k)}\ind_{\theta_2\leq p^n_{s_-}(j,\ell|k-1)}
\nonumber\\[-8pt]\\[-8pt]
&&\hspace*{79.3pt}{} \times\ind_{\theta
_3\leq\rho(u|j+1,n\muISn_{s_-})}
\ind_{\theta_4\leq\rho(v|\ell,n\muRSn_{s_-})}
\,\rmd {Q}^1\nonumber\\
&&{} -\frac{1}{n} \int_0^t\int_{\N} \delta_{\gamma
_i(n\muISn_{s_-})} \ind_{i\in{\II}^n_{s_-}}\,\rmd {Q}^2,\nonumber
\\
\label{eqmuRrenorm}
\muRSn_t &=& \muRSn_0+\frac{1}{n}
\int_0^t\int_{E_1} \Biggl(\sum_{i=1}^{R^n_{s_-}} \bigl(\delta_{\gamma
_{i}(n\muRSn_{s_-})-v_{i}}-\delta_{\gamma_{i}(n\muRSn_{s_-})}\bigr)\Biggr) \nonumber\\
&&\hspace*{81pt}{} \times\ind_{\theta_1\leq\lambda^n_{s_-}(k)n\muSn
_{s_-}(k)}\ind_{\theta_2\leq
p^n_{s_-}(j,\ell|k-1)}\nonumber\\[-8pt]\\[-8pt]
&&\hspace*{81pt}{} \times\ind_{\theta_3\leq\rho(u|j+1,n\muISn_{s_-})}
\ind_{\theta_4\leq\rho(v|\ell,n\muRSn_{s_-})}\,\rmd {Q}^1
\nonumber\\
&&{} + \frac{1}{n} \int_0^t\int_{\N} \delta_{\gamma
_i(n\muISn_{s_-})} \ind_{i\in{\II}^n_{s_-}}\,\rmd {Q}^2,\nonumber
\end{eqnarray}
where we denote for all $s\ge0$,
%
%
\begin{eqnarray}
\label{eqdeflambdan}
\lambda^n_s(k)&=&rk\frac{\NISnn_s}{\NSnn_s}\quad\mbox{and}\nonumber\\[-8pt]\\[-8pt]
p^n_s(j,\ell\mid k-1)&=&\frac{{\NISnn_{s}-1 \choose j} {\NRSnn_{s}
\choose\ell} {\NSnn_{s}-\NISnn_{s}-\NRSnn_{s}\choose k-1-j-\ell
}}{{\NSnn_{s}-1 \choose k-1}}.\nonumber
\end{eqnarray}

Several semimartingale decompositions will be useful in the
sequel. We focus on $\muISn$ but similar decompositions hold
for $\muSn$ and $\muRSn$, which we do not detail since
they can be deduced by direct adaptation of the following computation.
\begin{prop}
For all $f \in\mathcal B_b(\N)$, for all $t\ge0$,
%
%
\begin{equation}\label{musif}
\bigl\langle\muISn_t,f\bigr\rangle=\sum_{k\in\N}f(k)
\muISn_0(k)+A^{(n),{\II}
\sscS,f}_t+M^{(n),{\II}\sscS,f}_t,
\end{equation}
where the finite variation part $A^{(n),{\II}\sscS,f}_t$ of $\langle
\mu
^{(n),{\II}\sscS}_t,f\rangle$ reads
%
%
\begin{eqnarray}
\label{defAnSI}\qquad
A^{(n),{\II}\sscS,f}_t
&=& \int_0^t\sum_{k\in\N}\lambda
^n_{s}(k)\muSn_s(k)
\sum_{j+\ell+1\leq k}p_s^n(j,\ell|k-1)\nonumber\\
&&\hspace*{15pt}\hspace*{13.3pt}{}\times \sum_{u\in\U} \rho
(u|j+1,\muISnn_{s})\nonumber\\
&&\hspace*{15pt}\hspace*{13.3pt}{}\times \Biggl(f\bigl(k-(j+1+\ell)\bigr)
\\
&&\hspace*{46pt}{}+
\sum_{i=1}^{I^n_s} \bigl(f\bigl(\gamma_{i}(\muISnn_{s}) - u_i\bigr)-f(\gamma
_{i}(\muISnn_{s}))\bigr)\Biggr) \,\rmd {s}\nonumber\\
&&{}- \int_0^t \beta\bigl\langle \muISn_s,f
\bigr\rangle \,\rmd s,\nonumber
\end{eqnarray}
and where the martingale part $M^{(n),{\II}\sscS,f}_t$ of
$\langle\mu^{(n),{\II}\sscS}_t,f\rangle$ is a square integrable
martingale starting from 0 with quadratic variation
\begin{eqnarray*}
\bigl\langle
M^{(n),{\II}\sscS,f}\bigr\rangle_t 
&=& \frac{1}{n}\int_0^t \beta \bigl\langle \muISn_s,f^2
\bigr\rangle \,\rmd s\\
&&{}+\frac{1}{n} \int_0^t\sum_{k\in\N}
\lambda^n_s(k)\muSn_s(k)\sum_{j+\ell+1\leq
k}p^n_s(j,\ell|k-1) \\
&&\qquad\hspace*{29.4pt}{}\times\sum_{u\in\mathcal{U}}
\rho(u|j+1,\muISnn_{s})\\
&&\qquad\hspace*{29.4pt}{} \times\Biggl(f\bigl(k-(j+1+\ell)\bigr)\\
&&\hspace*{69.4pt}{}+\sum_{i=1}^{I^n_s}\bigl(f\bigl(\gamma_{i}(\muISnn
_{s}) - u_i\bigr)-f(\gamma_{i}(\muISnn_{s}))\bigr)\Biggr)^2 \,\rmd {s}.%
\end{eqnarray*}
\end{prop}
\begin{pf}
The proof proceeds from (\ref{eqmuIrenorm}) and standard stochastic
calculus for jump processes (see, e.g., \cite{fourniermeleard}).
\end{pf}

\subsection{Convergence of the normalized process}

We aim to study the limit of the system when $n\rightarrow+\infty$.
We introduce the associated measure spaces. For any $\varepsilon\geq
0$ and $A>0$, we define the following closed set of
$\M_F(\N)$ as
%
%
\begin{equation}\label{eqdefM}
\mathcal{M}_{\varepsilon, A}=\{\nu\in\mathcal{M}_F(\N) ;
\langle\nu,\car+\chi^5\rangle\leq A\mbox{ and } \langle\nu,
\chi\rangle\ge\varepsilon\}
\end{equation}
%
and $\M_{0+,A}=\bigcup_{\varepsilon>0}\M_{\varepsilon,A}$.
Topological properties of these spaces are given in the \hyperref
[appendiceA]{Appendix}.\vadjust{\goodbreak}

In the proof, we will see that the epidemic remains large provided the
number of edges from $\cI$ to $\cS$ remains of the order of $n$.
Let us thus define, for all $\varepsilon>0$, $\varepsilon'>0$ and
$n\in\N^*$,
%
%
\begin{equation}
\label{eqdeftepsilon}
t_{\varepsilon'}:= \inf\{t \ge0, \langle\bmuIS_t,\chi \rangle <
\varepsilon' \}
\end{equation}
and
%
%
\begin{equation}
\label{eqdeftauepsilonn}
\tau^n_\varepsilon=\inf\bigl\{t\geq0, \bigl\langle\muISn_t,\chi \bigr\rangle <
\varepsilon\bigr\}.
\end{equation}

Our main result is the following theorem.
\begin{Thm}\label{propconvergencemunS}
Suppose that Assumption \ref{hypconvcondinit} holds with
%
%
\begin{equation}
\label{eqhypmoments}
\qquad
\bigl(\mu^{(n),\sscS}_0,\mu^{(n),\scIS}_0,\mu^{(n),\cRSr}_0\bigr) \mbox{ in }
(\M_{0, A})^3\qquad\mbox{for any }n\mbox{, with }\langle
\bmuIS_0,\chi\rangle>0.
\end{equation}
Then, we have:\vspace*{8pt}

(1) 
there exists a unique solution $(\bmuS,\bmuIS,\bmuRS)$ to the
deterministic system of measure-valued equations
(\ref{limitereseauinfiniS})--(\ref{limitereseauinfiniSR}) in
$\Co(\R_+,\M_{0,A}\times \M_{0+,A}\times \M_{0,A})$,

(2) 
when $n$ converges to infinity, the sequence
$(\muSn,\muISn,\muRSn)_{n\in \N }$ converges in distribution in
$\D(\R_+,\M_{0,A}^3)$ to $(\bmuS,\bmuIS,\bmuRS)$.
\end{Thm}


\begin{pf}
Uniqueness of the solution to
(\ref{limitereseauinfiniS})--(\ref{limitereseauinfiniSR}) is proved in
Step 2. For the proof of (2), since $\lim_{\varepsilon'\rightarrow
0}t_{\varepsilon'}=+\infty$,  it is sufficient to prove the results on
$\D([0,t_{\varepsilon'}],\M_{0,A}^3)$ for $\varepsilon'$ small enough.
In the sequel, we choose $0<\varepsilon<\varepsilon'<\langle
\bmuIS_0,\chi\rangle$.

\textit{Step} 1. Let us prove that $(\muSn,\muISn,
\muRSn)_{n\in\N^*}$ is tight. Let $t\in\R_+$ and $n\in\N^*$. By
hypothesis, we have that
%
%
\begin{eqnarray}\label{domination0}
&&\bigl\langle\muSn_t,\car+ \chi^5 \bigr\rangle+\bigl\langle
\muISn_t,\car+ \chi^5 \bigr\rangle+ \bigl\langle\muRSn_t,\car+ \chi^5
\bigr\rangle\nonumber\\[-8pt]\\[-8pt]
&&\qquad\leq\bigl\langle\muSn_0,\car+ \chi^5 \bigr\rangle+\bigl\langle
\muISn_0,\car+ \chi^5\bigr\rangle\le 2A.
\nonumber
\end{eqnarray}
Thus, the sequences $(\muSn_t)_{n\in\N^*}$, $(\muISn_t)_{n\in\N
^*}$ and $(\muRSn_t)_{n\in\N^*}$ are tight for each $t\in\R_+$.
Now by the criterion of Roelly \cite{roelly}, it remains to prove that
for each $f\in\mathcal{B}_b(\N)$, the sequence $(\langle\muSn
_{\cdot},f\rangle, \langle\muISn_{\cdot},f\rangle, \langle\muRSn
_{\cdot},f\rangle
)_{n\in\N^*} $ is tight in $\D(\R_+,\R^3)$. Since we have
semimartingale decompositions of these processes, it is sufficient, by
using the Rebolledo criterion, to prove that the finite variation part
and the bracket of the martingale satisfy the Aldous criterion (see,
e.g., \cite{joffemetivier}). We only prove that $\langle\muISn
_{\cdot},f\rangle$ is tight. For the other components, the
computations are similar.

The Rebolledo--Aldous criterion is satisfied if for all
$\alpha>0$ and $\eta>0$ there exists $n_0\in\N$ and $\delta>0$
such that for all $n>n_0$ and for all stopping times~$S_n$ and $T_n$
such that $S_n<T_n<S_n+\delta$,
%
%
\begin{eqnarray}\label{criterealdous}
\PP\bigl(\bigl|A^{(n),{\II}\sscS,f}_{T_n}-A^{(n),{\II}\sscS,f}_{S_n}\bigr|>\eta
\bigr)&\leq&
\alpha\quad \mbox{and }\nonumber\\[-8pt]\\[-8pt]
\PP\bigl(\bigl|\bigl\langle
M^{(n),{\II}\sscS,f}\bigr\rangle_{T_n}-\bigl\langle
M^{(n),{\II}\sscS,f}\bigr\rangle_{S_n}\bigr|>\eta\bigr)&\leq&
\alpha.\nonumber
\end{eqnarray}

For the finite variation part,
\begin{eqnarray*}
&&\E\bigl[\bigl|A^{(n),{\II}\sscS,f}_{T_n}-A^{(n),{\II}\sscS,f}_{S_n}\bigr|\bigr]\\
&&\qquad\leq
\E\biggl[\int_{S_n}^{T_n} \beta \|f\|_\infty \bigl\langle
\muSIn_s,1\bigr\rangle \,\rmd s\biggr]\\
&&\qquad\quad{} +
\E\biggl[\int_{S_n}^{T_n} \sum_{k\in\N} \lambda^n_s(k) \muSn_s(k)
\sum_{j+\ell\leq k-1} p_s^n(j,\ell|k-1) (2j+1)\|f\|_\infty\,\rmd {s}\biggr].
\end{eqnarray*}
The term $\sum_{j+\ell\leq k-1}j p^n_s(j,\ell|k-1)$ is the mean
number of links to $\I^n_{s_-}$ that the newly infected individual
has, given that this individual is of degree $k$.
It is bounded by
$k$. Then, with (\ref{eqdeflambdan}),
\begin{eqnarray*}
&&
\E\bigl[\bigl|A^{(n),{\II}\sscS,f}_{T_n}-A^{(n),{\II}\sscS,f}_{S_n}\bigr|\bigr]
\\
&&\qquad \leq
\delta \E\bigl[\beta \|f\|_\infty \bigl(S^{(n)}_0+I^{(n)}_0\bigr)
+ r \|f\|_\infty \bigl\langle
\mu_0^{(n),\sscS},2\chi^2+3\chi\bigr\rangle\bigr],
\end{eqnarray*}
by using that the number of infectives is bounded by the size of the
population and that $\mu_s^{(n),\sscS}(k)\le\muSn_0(k)$ for all $k$ and
$s\geq0$. From (\ref{eqhypmoments}), the right-hand side is finite.
Using Markov's inequality,
\[
\PP\bigl(\bigl|A^{(n),{\II}\sscS,f}_{T_n}-A^{(n),{\II}\sscS,f}_{S_n}\bigr|>\eta
\bigr)\leq
\frac{(5 r+2\beta)A\delta
\|f\|_\infty}{\eta},
\]
which is smaller than $\alpha$ for $\delta$
small enough.

We use the same arguments for the bracket of the martingale
%
%
\begin{eqnarray}
\label{etape3}
&&\E\bigl[\bigl|\bigl\langle M^{(n),{\II}\sscS,f}\bigr\rangle_{T_n}-\bigl\langle M^{(n),{\II
}\sscS
,f}\bigr\rangle_{S_n}\bigr|\bigr]\nonumber\\
&&\qquad
\le \E\biggl[\frac{\delta \beta \|f\|^2_\infty (S^{(n)}_0+I^{(n)}_0)}{n}
+ \frac{\delta r\|f\|_\infty^2 \langle \muSn_0,\chi(2\chi+3)^2\rangle}{n}\biggr]\\
&&\qquad
\le \frac{(25 r+2\beta)A\delta \|f\|_\infty^2}{n},\nonumber
\end{eqnarray}
using assumption (\ref{eqhypmoments}). The right-hand side can be made
smaller than $\eta\alpha$ for a small enough $\delta$, so the
second inequality of (\ref{criterealdous}) follows again from Markov's
inequality. By \cite{roelly}, this provides the
tightness in $\D(\R_+,\mathcal{M}_{0,A}^3)$.

By Prohorov theorem (e.g., \cite{ethierkurtz}, page 104) and Step 1, the
distributions of $(\muSn,\muISn, \muRSn)$, for $n\in\N^*$, form a
relatively compact family of bounded measures on $\D(\R_+,\M_{0,A}^3)$,
and so do the laws of the stopped processes $(\muSn_{\cdot\wedge
\tau^n_\varepsilon},\muISn_{\cdot\wedge\tau^n_\varepsilon},
\muRSn_{\cdot\wedge\tau^n_\varepsilon})_{n\in\N^*}$ [recall
(\ref
{eqdeftauepsilonn})]. Let
$\bar{\mu}:=(\bmuS,\bmuIS,\bmuRS)$ be a~limiting point in
$\Co(\R_+,\mathcal{M}_{0,A}^3)$
of the sequence of stopped processes
and let us consider a subsequence again
denoted by
$\mu^{(n)}:=(\muSn,\muISn, \muRSn)_{n\in\N^*}$, with an abuse
of\vadjust{\goodbreak}
notation, and that converges to
$\bar\mu$. Because the limiting values are continuous, the
convergence of $(\mu^{(n)})_{n\in\N^*}$ to $\bar\mu$ holds for the
uniform convergence on every compact subset of $\R_+$ 
(e.g., \cite{billingsley}, page 112).

Now, let us define for all $t\in\R_+$ and for all bounded function
$f$ on $\N$, the mappings $\Psi_t^{\sscS,f}$, $\Psi_t^{{\II}\sscS
,f}$ and
$\Psi_t^{\sRr\sscS,f}$ from $\D(\R_+,\mathcal{M}_{0,A}^3)$ into
$\D(\R_+,\R)$ such that~(\ref{limitereseauinfiniS})--(\ref{limitereseauinfiniSR}) read
%
%
\begin{eqnarray}
\label{eqdefPsi}\qquad
&&(\langle\bmuS_t,f \rangle,\langle\bmuIS_t,f \rangle,\langle
\bmuRS_t,f \rangle)
\nonumber\\[-8pt]\\[-8pt]
&&\qquad=(\Psi^{\sscS,f}_t(\bmuS,\bmuIS,\bmuRS),\Psi^{{\II
}\sscS
,f}_t(\bmuS
,\bmuIS,\bmuRS),\Psi^{\sRr\sscS,f}_t(\bmuS,\bmuIS,\bmuRS)).
\nonumber
\end{eqnarray}
Our purpose is to prove that the limiting values are the unique
solution of equations
(\ref{limitereseauinfiniS})--(\ref{limitereseauinfiniSR}).

Before proceeding to the proof, a remark is in order. A natural
way of reasoning would be to prove that $\Psi^{\sscS,f},
\Psi^{{\II}\sscS,f}$ and $\Psi^{\sRr\sscS,f}$ are Lipschitz
continuous in some
spaces of measures. It turns out that this can only be done by considering
the set of measures with moments of any order, which is a~set too small for applications. We circumvent this difficulty by first proving
that the mass and the first two moments of any solutions of the system
are the same. Then, we prove that the generating
functions of these measures satisfy a partial differential equation
known to have a unique solution.

\textit{Step} 2. 
We now prove that the differential system
(\ref{limitereseauinfiniS})--(\ref{limitereseauinfiniSR}) has at most
one solution in ${\mathcal C}(\R_+, \M_{0,A}\times \M_{0+,A}\times
\M_{0,A})$. It is enough to prove the result in ${\mathcal C}([0,T],
\M_{0,A}\times \M_{\varepsilon,A}\times \M_{0,A})$ for all
$\varepsilon>0$ and $T>0$. Let
$\bar{\mu}^i=(\bar{\mu}^{\sscS,i},\bar{\mu}^{{\II}\sscS,i},\bar{\mu}^{
\sRr\sscS,
i})$,
$i\in \{1,2\}$ be two solutions of
(\ref{limitereseauinfiniS})--(\ref{limitereseauinfiniSR}) in this
space, started with the same initial conditions. Set
\begin{eqnarray*}
\Upsilon_t &=& \sum_{j=0}^3 |\langle\bar{\mu}^{\sscS,1}_t,\chi
^j\rangle-\langle
\bar{\mu}^{\sscS,2}_t,\chi^j\rangle| \\
&&{} + \sum_{j=0}^2(|\langle
\bar{\mu}^{\scIS,1}_t,\chi^j\rangle-\langle
\bar{\mu}^{\scIS,2}_t,\chi^j\rangle|+|\langle
\bar{\mu}^{\cRSr,1}_t,\chi^j\rangle-\langle
\bar{\mu}^{\cRSr,2}_t,\chi^j\rangle|
).
\end{eqnarray*}

Let us first remark that for all $0\leq t <T$,
$\bNS_t\geq \bNIS_t>\varepsilon$ and then
%
%
\begin{eqnarray}\label{etape81}\quad
|\bar{p}^{\II,1}_t-\bar{p}^{\II,2}_t|&=&
\biggl|\frac{\bar{N}^{\scIS,1}_t}{\bar{N}^{\sscS,1}_t}-\frac{\bar
{N}^{\scIS,2}_t}{\bar{N}^{\sscS,2}_t}\biggr|\nonumber\\
&\leq&\frac{A}{\varepsilon^2}|\bar{N}^{\sscS,1}_t-\bar{N}^{\sscS,2}_t|+
\frac{1}{\varepsilon}|\bar{N}^{\scIS,1}_t-\bar{N}^{\scIS,2}_t|\nonumber\\[-8pt]\\[-8pt]
&=& \frac{A}{\varepsilon^2}|\langle\bar{\mu}_t^{\sscS,1}, \chi
\rangle-\langle\bar{\mu}_t^{\sscS,2}, \chi \rangle|+\frac
{1}{\varepsilon}|\langle\bar{\mu}_t^{\scIS,1}, \chi \rangle
-\langle\bar{\mu}_t^{\scIS,2}, \chi \rangle|\nonumber\\
&\le&\frac
{A}{\varepsilon
^2}\Upsilon_t.
\nonumber
\end{eqnarray}
The same computations show a similar result for\vadjust{\goodbreak} $
|\bar{p}^{\sscS,1}_t-\bar{p}^{\sscS,2}_t|$. 

Using that $\bar{\mu}^i$ are solutions to
(\ref{limitereseauinfiniS})--(\ref{limitereseauinfiniSI}) let us show
that $\Upsilon$ satisfies a~Gronwall inequality which implies
that it is equal to 0 for all $t\le T$.
For the degree distributions of the susceptible individuals, we have
for $p\in\{0,1,2,3\}$ and $f=\chi^p$ in (\ref{limitereseauinfiniS})
\begin{eqnarray*}
|\langle\bar{\mu}^{\sscS,1}_t,\chi^p\rangle-\langle\bar{\mu
}^{\sscS,2}_t,\chi^p\rangle| &= & \biggl|\sum_{k\in\N} \bar{\mu
}_0^{\sscS
}(k) k^p (e^{-r\int_0^t \bar{p}^{\II,1}_s \,\rmd s}-e^{-r\int_0^t \bar
{p}^{\II,2}_s \,\rmd s})\biggr|\nonumber\\
&\leq& r \sum_{k\in\N} k^p \bar{\mu}_0^{\sscS}(k) \int_0^t
|\bar
{p}^{\II,1}_s-\bar{p}^{\II,2}_s |\,\rmd s\\
&\leq& r\frac{A^2}{\varepsilon^2}\int_0^t \Upsilon_s \,\rmd s
\end{eqnarray*}
by using (\ref{etape81}) and the fact that $\bmuS_0\in\mathcal
{M}_{0,A}$.

For $\bmuSI$ and $\bmuSR$, we use (\ref{limitereseauinfiniSI}) and
(\ref{limitereseauinfiniSR}) with the
functions $f=\chi^0$, $f=\chi$ and $f=\chi^2$. We proceed here with
only one of the
computations, others can be done similarly. From (\ref{limitereseauinfiniSI}),
\begin{eqnarray*}
&&
\langle\bar{\mu}_t^{\scIS,1}, \mathbf{1} \rangle-\langle\bar
{\mu}_t^{\scIS,2}, \mathbf{1} \rangle\\[3pt]
&&\qquad=\beta\int_0^t\langle\bar
{\mu}_s^{{\II}\sscS,1}-\bar{\mu}_s^{\scIS,2}, \mathbf{1} \rangle
\,\rmd {s}
+r \int_0^t (\bar{p}_s^{\II,1}\langle\bar{\mu}^{\sscS,1}_s,
\chi \rangle-\bar{p}_s^{\II,2}\langle\bar{\mu}^{\sscS,2}_s,
\chi \rangle)\,\rmd {s}.
\end{eqnarray*}
Hence, with (\ref{etape81}),
%
\[
| \langle\bar{\mu}_t^{\scIS,1}-\bar{\mu}_t^{\scIS,2}, \mathbf
{1} \rangle|\le
C(\beta,r,A,\varepsilon) \int_0^t \Upsilon_s\,\rmd {s}.
\]
By analogous computations for the other quantities, we then show that
\[
\Upsilon_t\le C'(\beta,r,A,\varepsilon) \int_0^t \Upsilon_s\,\rmd {s},
\]
hence, $\Upsilon\equiv0$.
It follows that for all
$t<T$, and for all $j\in
\{0,1,2\}$,
%
%
\begin{equation}\label{etape10}
\langle\bar{\mu}^{\sscS,1}_t,\chi^j\rangle=\langle
\bar{\mu}^{\sscS,2}_t,\chi^j\rangle\quad\mbox{and}\quad \langle
\bar{\mu}^{{\II}\sscS,1}_t,\chi^j\rangle=\langle
\bar{\mu}^{{\II}\sscS,2}_t,\chi^j\rangle,
\end{equation}
and in particular, $\bar{N}^{\sscS,1}_t=\bar{N}^{\sscS,2}_t$ and
$\bar{N}^{{\II}\sscS,1}_t=\bar{N}^{{\II}\sscS,2}_t$. This implies
that $
\bar{p}^{\sscS,1}_t=\bar{p}^{\sscS,2}_t$, $\bar{p}^{{\II
},1}_t=\bar
{p}^{{\II},2}_t$ and $\bar{p}^{\sRr,1}_t=\bar{p}^{\sRr,2}_t$. From
(\ref{limitereseauinfiniS}) and the continuity of the solutions to
(\ref{limitereseauinfiniS})--(\ref{limitereseauinfiniSR}), pathwise
uniqueness holds for $\bmuS$ a.s.

Our purpose is now to prove that
$\bar{\mu}^{\scIS,1}=\bar{\mu}^{\scIS,2}$. Let us introduce the
following generating functions: for any
$t\in
\R_+$, $i\in\{1,2\}$ and $\eta\in[0,1)$,
\[
\mathcal{G}^i_t(\eta)=\sum_{k\geq0}\eta^k
\bar{\mu}^{{\II}\sscS,i}_t(k).
\]
Since we already know these
measures do have the same total mass, it boils down to prove that $
\mathcal{G}^1\equiv\mathcal{G}^2$.
Let us define
%
%
\begin{eqnarray}\label{etape11}\qquad
H(t,\eta) &=& \int_0^t \sum_{k\in\N}rk\bpI_s\mathop{\sum_{j,
\ell, m\in
\N}}_{j+\ell+m=k-1}\pmatrix{k-1 \cr j,\ell,m}(\bpI_s)^j (\bpR
_s)^\ell
(\bpS_s)^m \eta^m \bmuS_s(k) \,\rmd s,\nonumber\\[-8pt]\\[-8pt]
K_t&=&\sum_{k\in
\N}rk\bpI_t(k-1)\bpR_t\frac{\bmuS_t(k)}{\bNIS_t}.\nonumber
\end{eqnarray}
The latter quantities are, respectively, of class $\Co^1$ and
$\Co^0$ with respect to time $t$ and are well defined and bounded on
$[0,T]$. Moreover, $H$ and $K$ do not depend
on the chosen solution because of (\ref{etape10}). Applying
(\ref{limitereseauinfiniSI}) to $f(k)=\eta^k$ yields
\begin{eqnarray*}
\mathcal{G}^i_t(\eta) &=& \mathcal{G}^i_0(\eta)+H(t,\eta)+ \int
_0^t \biggl(K_s \sum_{k'\in\N^*}(\eta^{k'-1}-\eta^{k'})k'
\bar{\mu}^{\scIS,i}
_s(k')-\beta\mathcal{G}^i_s(\eta)\biggr) \,\rmd s\\
&=& \mathcal{G}^i_0(\eta)+H(t,\eta)+ \int_0^t \bigl(K_s
(1-\eta)\partial_\eta\mathcal{G}^i_s(\eta) -\beta
\mathcal{G}^i_s(\eta)\bigr) \,\rmd s.
\end{eqnarray*}
Then, the
functions $t\mapsto\widetilde{\mathcal{G}}^i_t(\eta)$ defined by
$\widetilde{\mathcal{G}}^i_t(\eta)=e^{\beta t}\mathcal{G}^i_t(\eta)$,
$i\in\{1,2\}$, are solutions of the following transport equation:
%
%
\begin{equation}
\partial_t g(t,\eta)-(1-\eta)K_t \partial_\eta
g(t,\eta)=\partial_t H(t,\eta)e^{\beta t}.
\end{equation}
In view of the regularity of $H$ and $K$, it is known that this
equation admits a unique solution (see, e.g., \cite{evans}). Hence,
$\mathcal{G}^1_t(\eta)=\mathcal{G}^2_t(\eta)$ for all $t\in\R_+$ and
$\eta\in[0,1)$. The same method applies to $\bmuRS$. Thus, there is
at most
one solution to the differential system (\ref
{limitereseauinfiniS})--(\ref{limitereseauinfiniSR}).

\textit{Step} 3. We now show that $\mu^{(n)}$ nearly satisfies
(\ref{limitereseauinfiniS})--(\ref{limitereseauinfiniSR}) as $n$ gets
large. Recall (\ref{musif}) for a bounded function
$f$ on $\N$. To identify the limiting values, we establish that for
all $n\in\N^*$ and all $t\ge0$,
%
%
\begin{equation}\label{decompo}
\bigl\langle\muISn_{t\wedge\tau^n_\varepsilon},f\bigr\rangle=\Psi_{t\wedge
\tau^n_\varepsilon}^{{\II}\sscS,f}\bigl(\mu^{(n)}\bigr)+\Delta_{t\wedge
\tau^n_\varepsilon}^{n,f}+M^{(n),{\II}\sscS,f}_{t\wedge
\tau^n_\varepsilon},
\end{equation}
where $M^{(n),{\II}\sscS,f}$ is defined in (\ref{musif}) and where
$\Delta
_{\cdot\wedge\tau^n_\varepsilon}^{n,f}$ converges to 0 when
$n\rightarrow+\infty$, in probability and uniformly in $t$ on compact
time intervals.

Let us fix $t\in\R_+$. Computations similar to (\ref{etape3}) give
%
%
\begin{equation}\label{etape777}
\E\bigl(\bigl(M^{(n),{\II}\sscS,f}_t\bigr)^2\bigr)=\E\bigl(\bigl\langle
M^{(n),{\II}\sscS,f}\bigr\rangle_t\bigr) \leq
\frac{(25r+2\beta) A t
\|f\|_\infty^2}{n}.
\end{equation}
Hence, the sequence $(M^{(n),{\II}\sscS,f}_t)_{n\in\N}$ converges in
$L^2$ and in probability to zero (and in $L^1$ by Cauchy--Schwarz
inequality).

We now consider the finite variation part of (\ref{musif}), given
in (\ref{defAnSI}). The sum in (\ref{defAnSI}) corresponds to the
links to $\I$ that the new infected individual has. We separate this
sum into cases where the new infected individual only has simple edges
to other
individuals of $\cI$, and cases where multiple edges
exist. The latter term is expected to vanish for large populations:
%
%
\begin{equation}\label{eqshadow1}
A^{(n),{\II}\sscS,f}_t
%
= B^{(n),{\II}\sscS,f}_t + C^{(n),{\II}\sscS,f}_t,
\end{equation}
where
%
%
\begin{eqnarray}\label{termeBn}
B^{(n),{\II}\sscS,f}_t &=&
-\int_0^t \beta\bigl\langle \muSIn_s,f\bigr\rangle\,\rmd s\nonumber\\
&&{}+
\int_0^t \sum_{k\in\N}
\lambda^n_s(k)\muSn_s(k)
\sum_{j+\ell+1\leq k}p^n_s(j,\ell
|k-1)\nonumber\\[-8pt]\\[-8pt]
&&{}
\times\Biggl\{f\bigl(k-(j+1+\ell)\bigr)
+ \mathop{\sum_{u\in\mathcal{U}(j+1,\muISnn_s);}}_{\forall
i\leq I^n_{s_-}, u_i\leq1
}\rho(u|j+1,\muISnn_s)\nonumber\\
&&\hspace*{77pt}{}\times
\sum_{i=0}^{I^n_{s_-}}\bigl(f\bigl(\gamma_{i}(\muSInn_{s_-})
-
u_i\bigr)-f(\gamma_{i}(\muSInn_{s_-}))\bigr)\Biggr\} \,\rmd {s}
\nonumber
\end{eqnarray}
and
%
%
\begin{eqnarray}\label{termeCn}
C^{(n),{\II}\sscS,f}_t &=& \int_0^t \sum_{k\in\N} \lambda
^n_s(k)\muSn
_s(k)
\sum_{j+\ell+1\leq
k}p^n_s(j,\ell|k-1)\nonumber\\
&&\hspace*{18pt}\hspace*{11pt}{}\times\mathop{\sum_{u\in\mathcal{U}(j+1,\muISnn_s);}}_{
\exists i\leq I^n_{s_-}, u_i>1}\rho(u|j+1,\muISnn_s)\\
&&\hspace*{18pt}\hspace*{11pt}{}\times\sum_{i=1}^{I^n_{s_-}}\bigl(f\bigl(\gamma_{i}(\muSInn_{s_-}) -
u_i\bigr)-f(\gamma_{i}(\muSInn_{s_-}))\bigr) \,\rmd {s}.
\nonumber
\end{eqnarray}
%
We first show that $C^{(n),\sscS{\II},f}_t$ is a negligible term. Let
$q^n_{j,\ell,s}$ denote the probability that the newly infected
individual at time $s$ has a double (or of higher order) edge to some
alter in $\I_{s_-}^n$, given $j$ and $\ell$. The probability to have
a~multiple edge to a given infectious $i$ is less than the number of
couples of edges linking the newly infected to $i$, times the
probability that these two particular edges linking $i$ to a
susceptible alter at $s_-$ actually lead to the newly infected. Hence,
%
%
\begin{eqnarray}\label{estimatep2}
q^n_{j,\ell,s} &=& \mathop{\sum_{u\in\mathcal{U}(j+1,\muISnn
_s);}}_{\exists i\leq I^n_{s_-}, u_i>1}\rho(u|j+1,\muISnn_s)\nonumber
\\
&\leq& \pmatrix{j\cr2} \sum_{x\in
{\II}^n_{s_-}}\frac{d_x(\sscS^n_{s_-})(d_x(\sscS
^n_{s_-})-1)}{\NISnn
_{s_-}(\NISnn_{s_-}-1)}\nonumber\\[-8pt]\\[-8pt]
&=& \pmatrix{j \cr2}\frac{1}{n}\frac{\langle\muISn_{s_-},\chi(\chi
-1)\rangle}{\NISn_{s_-}(\NISn_{s_-}-1/n)}\nonumber\\
&\leq& \pmatrix{j\cr2}
\frac{1}{n}\frac{A}{\varepsilon(\varepsilon-1/n)} \qquad\mbox{if
}s<\tau^n_\varepsilon\mbox{ and }n>1/\varepsilon.\nonumber
\end{eqnarray}
Then, since for all $u\in\mathcal{U}(j+1,\muISnn_s)$,
%
%
\begin{equation}\label{ancienneeq325}
\Biggl|\sum_{i=1}^{I^n_{s_-}}\bigl(f\bigl(\gamma_{i}(\muSInn_{s_-}) - u_i\bigr)-f(\gamma
_{i}(\muSInn_{s_-}))\bigr) \Biggr|\leq2(j+1)\|f\|_\infty,
\end{equation}
we have by (\ref{estimatep2}) and (\ref{ancienneeq325}), for
$n>1/\varepsilon$,
%
%
\begin{eqnarray}
\label{etape5}\hspace*{32pt}
\bigl|C^{(n),{\II}\sscS,f}_{t\wedge
\tau^n_\varepsilon}\bigr|
&\leq& \int_0^{t\wedge\tau^n_\varepsilon} \sum_{k\in\N} rk \muSn
_s(k)\sum_{j+\ell+1\leq k}p^n_s(j,\ell|k-1)
2(j+1)\|f\|_\infty\nonumber\\
&&\hspace*{132.5pt}{}\times\frac
{j(j-1)A}{2n\varepsilon(\varepsilon-1/n)} \,\rmd {s}\\
&\leq& \frac{A rt\|f\|_\infty}{n \varepsilon
(\varepsilon-1/n)}\bigl\langle\muSn_0,\chi^4\bigr\rangle,\nonumber
\end{eqnarray}
which tends to zero in view of (\ref{eqhypmoments})
and to the fact that $\muSn_s$ is dominated by $\mu_0^{(n),\sscS}$ for
all $s\geq0$ and $n\in\N^* $.

We now aim at proving that
$B^{(n),{\II}\sscS,f}_{\cdot\wedge\tau^n_\varepsilon}$ is somewhat
close to
$\Psi_{\cdot\wedge\tau^n_\varepsilon}^{{\II}\sscS,f}(\mu^{(n)})$.
First, notice
that
%
%
\begin{eqnarray}
\label{eqshadow2}
&&\mathop{\sum_{u\in\mathcal{U}(j+1,\muISnn_s);}}_{\forall i\leq I^n_{s_-},
u_i \le1}\rho(u|j+1,\muSInn_s)\nonumber\\
&&\quad\hspace*{48.7pt}{}\times\sum_{i=1}^{I_{s_-}^n}\bigl(f\bigl(\gamma
_{i}(\muSInn_{s_-})-u_i\bigr)-f(\gamma_{i}(\muSInn_{s_-}))\bigr)\nonumber\\
&&\qquad= \mathop{\sum_{u\in({\II}_{s_-}^n)^{j+1}}}_{u_0\not=
\cdots
\not=u_j} \biggl(\frac{\prod_{k=0}^j
d_{u_k}(\sscS_{s}^n)}{N_{s_-}^{n,{\II}\sscS}\cdots
(N_{s_-}^{n,\sscS{\II}}-(j+1))}\biggr)\nonumber\\
&&\hspace*{39pt}\qquad\quad{}\times
\sum_{m=0}^j\bigl(f\bigl(d_{u_m}(\sscS^n_{s_-})
- 1\bigr)-f(d_{u_m}(\scS_{s_-}^n))\bigr)\nonumber\\
&&\qquad= \sum_{m=0}^j\mathop{\sum_{u\in({\II}
_{s_-}^n)^{j+1}}}_{u_0\not=
\cdots\not=u_j} \biggl(\frac{\prod_{k=0}^j
d_{u_k}(\sscS_{s}^n)}{N_{s_-}^{n,{\II}\sscS}\cdots
(N_{s_-}^{n,\sscS{\II}}-(j+1))}\biggr)\\
&&\qquad\quad\hspace*{56.5pt}{}\times\bigl(f\bigl(d_{u_m}(\scS_{s_-}^n)
- 1\bigr)-f(d_{u_m}(\scS_{s_-}^n))\bigr)\nonumber\\
&&\qquad= \sum_{m=0}^j \biggl(\sum_{x\in{\II}^n_{s_-}}
\frac{d_{x}(\sscS^n_{s_-})}{\NISnn_{s_-}}\bigl(f\bigl(d_{x}(\scS^n_{s_-})
- 1\bigr)-f(d_{x}(\scS^n_{s_-}))\bigr)\biggr)\nonumber\\
&&\qquad\quad\hspace*{15.6pt}{}\times\biggl(\mathop{\sum_{u\in
({\II}_{s_-}^n\setminus\{x\})^{j}}}_{u_0\not= \cdots
\not=u_{j-1}}\frac{\prod_{k=0}^{j-1}
d_{u_k}(\sscS_{s}^n)}{(N_{s_-}^{n,{\II}\sscS}-1)\cdots
(N_{s_-}^{n,{\II}\sscS
}-(j+1))}\biggr)\nonumber\\
&&\qquad= (j+1) \frac{\langle
\muISn_{s_-},\chi(\tau_1f-f)\rangle}{\NISn_{s_-}}(1-q^n_{j-1,\ell,s}),
\nonumber
\end{eqnarray}
where we recall that $\tau_1f(k)=f(k-1)$ for every
function $f$ on $\N$ and $k\in\N$.
In the third equality, we split the term $u_m$ from the other terms
$(u_{m'})_{m'\not= m}$. The last sum in the right-hand side of this
equality is the probability of drawing $j$ different infectious
individuals that are not $u_m$ and that are all different, hence,
$1-q^n_{j-1,\ell,s}$.

Denote for $t>0$ and $n\in\N$,
\begin{eqnarray*}
p^{n,{\II}}_t&=&\frac{\langle\muISnn_t,\chi \rangle-1}{\langle
\mu _t^{n,\sscS },\chi \rangle-1},\\
p^{n,\sRr}_t&=&\frac{\langle\muRSnn_t,\chi \rangle}{\langle\mu
_t^{n,\sscS },\chi \rangle-1},\\
p^{n,\sscS}_t&=&\frac{\langle\muSnn_t,\chi \rangle-\langle\muISnn
_t,\chi \rangle-\langle\muRSnn_t,\chi \rangle}{\langle\mu
_t^{n,\sscS},\chi \rangle-1},
\end{eqnarray*}
the proportion of edges with infectious (resp., removed and
susceptible) alters and susceptible egos among all the edges with
susceptible egos but the contaminating edge. For all integers $j$ and
$\ell$ such that $j+\ell\leq k-1$ and $n\in\N^* $, denote by
\[
\tilde
p^n_t(j,\ell\mid
k-1)=\frac{(k-1)!}{j!(k-1-j-\ell)!\ell!}(p_t^{n,{\II
}})^{j}(p_t^{n,\sRr
})^{\ell}(p_t^{n,\sscS})^{k-1-j-\ell},
\]
the probability that the multinomial variable counting the number of
edges with infectious, removed and susceptible alters, among $k-1$
given edges, equals
$(j,\ell,k-1-j-\ell)$. We have that
%
%
\begin{equation}\label{etape778}
\bigl|\Psi_{t\wedge
\tau^n_\varepsilon}^{{\II}\sscS,f}\bigl(\mu^{(n)}\bigr)-B^{(n),{\II}\sscS
,f}_{t\wedge
\tau^n_\varepsilon}\bigr|
\leq \bigl|D^{(n),{\II}\sscS,f}_{t\wedge
\tau^n_\varepsilon}\bigr|+\bigl|E^{(n),{\II}\sscS,f}_{t\wedge
\tau^n_\varepsilon}\bigr|,
\end{equation}
where
\begin{eqnarray*}
D^{(n),{\II}\sscS,f}_{t} & = & \int_0^t \sum_{k\in\N} \lambda
^n_s(k)\muSn_s(k)\\
&&\hspace*{28.6pt}{}\times\sum_{j+\ell+1\leq k}\bigl(p^n_s(j,\ell|k-1)-\tilde
p^n_s(j,\ell|k-1)\bigr)\\
&&\hspace*{76pt}{} \times\biggl(f\bigl(k-(j+\ell+1)\bigr)\\
&&\hspace*{94.5pt}{} +
(j+1) \frac{\langle\muISn_{s_-},\chi(\tau_1f-f)\rangle}{\NISn
_{s_-}}\biggr) \,\rmd {s},\\
E^{(n),{\II}\sscS,f}_{t} & = & \int_0^t \sum_{k\in\N}
\lambda^n_s(k)\muSn_s(k)\\
&&\hspace*{28.6pt}{} \times\sum_{j+\ell+1\leq k}p^n_s(j,\ell|k-1)
(j+1)\\
&&\hspace*{79.2pt}\times\frac{\langle
\muISn_{s_-},\chi(\tau_1f-f)\rangle}{\NISn_{s_-}}q^n_{j-1,\ell,s}
\,\rmd {s}.
\end{eqnarray*}
First,
%
%
\begin{equation}\label{etape753}
\bigl|D^{(n),{\II}\sscS,f}_{t\wedge\tau^n_\varepsilon}\bigr|\leq\int
_0^{t\wedge
\tau^n_\varepsilon} \sum_{k\in\N} rk \alpha^n_s(k) \|f\|_\infty
\biggl(1+ \frac{2kA}{\varepsilon}\biggr) \muSn_s(k) \,\rmd {s},
\end{equation}
where for all $k\in\N$
\[
\alpha^n_t(k)=\sum_{j+\ell+1\leq
k} \bigl|p^n_t(j,\ell|k-1)-\tilde p^n_t(j,\ell|k-1) \bigr|.
\]
The
multinomial probability $\tilde p^n_s(j,\ell|k-1)$ approximates the
hypergeometric one, $p^n_s(j,\ell|k-1,s)$, as $n$ increases to
infinity, in view of the fact that the total population size,
$\langle\mu_0^{n,\sscS},\car \rangle+\langle\muISnn_0,\car
\rangle$, is of order $n$.
Hence, the right-hand side of (\ref{etape753}) vanishes by dominated
convergence.

On the other hand, using (\ref{estimatep2}),
%
%
\begin{eqnarray}\label{etape779}
\bigl|E^{(n),{\II}\sscS,f}_{t\wedge\tau^n_\varepsilon}\bigr| &\leq& \int
_0^{t\wedge\tau^n_\varepsilon} \sum_{k\in\N} rk^2 \muSn_s(k)
\frac{2\|f\|_\infty A}{\varepsilon}\frac{k^2 A}{2n\varepsilon
(\varepsilon-1/n)} \,\rmd {s}\nonumber\\[-8pt]\\[-8pt]
&\leq& \frac{A^3
rt\|f\|_\infty}{n\varepsilon^2(\varepsilon-1/n)},\nonumber
\end{eqnarray}
in view of (\ref{eqhypmoments}). Gathering
(\ref{etape777}), (\ref{eqshadow1}), (\ref{etape5}),
(\ref{etape778}), (\ref{etape753}) and (\ref{etape779}) concludes the
proof that
the rest of (\ref{decompo}) vanishes in probability uniformly over
compact intervals.

\textit{Step} 4. Recall\vspace*{-1pt} that in this proof, $\bar{\mu
}=(\bmuS,\bmuIS,\bmuRS)$ is the limit of $\mu^{(n)}_{\cdot\wedge \tau
^n_\varepsilon}=(\muSn_{\cdot\wedge\tau^n_\varepsilon},\muISn
_{\cdot\wedge\tau^n_\varepsilon},\muRSn_{\cdot\wedge\tau ^n_\varepsilon
})_{n\in\N^*}$, and recall that these processes take values in the
closed set $\M_{0,A}^3$. Our purpose is now to prove that $\bar{\mu }$
satisfy (\ref{limitereseauinfiniS})--(\ref{limitereseauinfiniSR}).
Using the Skorokhod representation theorem, there exists, on the same
probability space as~$\bar{\mu}$, a sequence,\vspace*{1pt} again denoted by
$(\mu^{(n)}_{\cdot\wedge \tau^n_\varepsilon})_{n\in \N^* }$
with\vspace*{1pt} an abuse of notation, with the same marginal
distributions as the original sequence, and that converges a.s. to~$\bar{\mu}$.\vspace*{1pt}

The maps $\nu_{\cdot}:=(\nu^1_{\cdot},\nu^2_{\cdot},\nu^3_{\cdot
})\mapsto\langle
\nu^{1}_{\cdot},\car\rangle/(\langle\nu_0^{1},\car\rangle
+\langle
\nu_0^{2},\car\rangle+\langle
\nu_0^{3},\car\rangle)$ [resp., $\langle\nu^{2}_{\cdot},\car
\rangle
/(\langle
\nu_0^{1},\car\rangle+\langle\nu_0^{2},\car\rangle+\langle
\nu_0^{3},\car\rangle)$ and
$\langle
\nu^{3}_{\cdot},\car\rangle/(\langle\nu_0^{1},\car\rangle
+\langle
\nu_0^{2},\car\rangle)+\langle
\nu_0^{3},\car\rangle$] are continuous from $\Co(\R_+,\mathcal
{M}_{0,A}\times
\M_{\varepsilon,A}\times\M_{0,A})$ into
$\Co(\R_+,\R)$.

Then, Lemma \ref{lemmaconti1} together with
the continuity of $(X^1_{\cdot},X^2_{\cdot})\mapsto X^1_{\cdot
}/X^2_{\cdot}$ from
$\C(\R_+,\R)\times\C(\R_+,\R^*)$ into
$\C(\R_+,\R)$ (see, e.g., \cite{Whitt80}), implies that the mapping
$\nu_{\cdot}\mapsto\langle\nu^1_{\cdot},\chi \rangle/\langle
\nu^2_{\cdot},\chi \rangle$ is continuous
from $\C(\R_+,\M_{0,A}\times\M_{\varepsilon,A}\times\M_{0,A})$ into
$\C(\R_+,\R)$.
The same argument yields the continuity of $\nu_{\cdot}\mapsto\ind
_{\langle\nu_{\cdot}^{1},\chi\rangle>\varepsilon}/\allowbreak\langle\nu
_{\cdot}^{2},\chi\rangle$ for the same spaces.

Lemma \ref{lemmaconti1} also provides the continuity of $\nu_{\cdot
}\mapsto
\langle\nu^2_{\cdot},\chi(\tau_1f-f) \rangle$ from $\C(\R_+,\M
_{0,A}\times\M
_{\varepsilon,A}\times\M_{0,A})$ into $\C(\R_+,\R)$ for bounded
function $f$ on $\N$.

Since, as well known, the mapping $y\in\D([0,t],\R)\mapsto\int_0^t
y_s \,\rmd {s}$ is continuous, we have proven the continuity of the
mapping $\Psi^f_t$ defined in (\ref{eqdefPsi}) on $\D(\R_+,\mathcal
{M}_{0,A}\times\mathcal{M}_{\varepsilon,A} \times\mathcal
{M}_{0,A})$.

By Lemma \ref{lemmaconti1} applied to $\varphi= \chi$, the process
$(\NISn_{\cdot\wedge\tau^n_\varepsilon})_{n\in\N^* }$ converges in
distribution to $\bNIS_{\cdot}=\langle\bmuIS_{\cdot},\chi\rangle
$. Since the
latter process is continuous, the convergence holds in
$(\D([0,T],\R_+),\|\cdot\|_{\infty})$ for any $T>0$ (see
\cite{billingsley}, page 112). As $y\in\D(\R_+,\R) \mapsto\inf
_{t\in
[0,T]} y(t)\in\R$ is continuous, we have a.s. that
\[
\inf_{t\in[0,T]} \bNIS_t = \lim_{n\rightarrow+\infty} \inf_{t\in
[0,T]} \NISn_{t\wedge\tau^n_\varepsilon} \qquad(\mbox{$\geq$}\varepsilon).
\]
%
We consider $\bar{t}_{\varepsilon'}=\inf\{t\in\R_+,
\bNIS_t\leq\varepsilon'\}$. A difficulty lies in the fact that we do not
yet know whether this time is deterministic. We have
a.s.
%
%
\begin{equation}
\label{etape2}
\varepsilon'\leq
\inf_{t\in[0,T]} \bNIS_{t\wedge \bar{t}_{\varepsilon'}} =
\lim_{n\rightarrow+\infty} \inf_{t\in[0,T]} \NISn_{t\wedge
\tau_\varepsilon^n \wedge \bar{t}_{\varepsilon'}}.
\end{equation}
Hence, using Fatou's lemma,
%
%
\begin{eqnarray}\label{etapelimitetaun}
1 &=& \PP\Bigl(\inf_{t\in[0,\bar t_{\varepsilon'}]}
\bNIS_{t}>\varepsilon\Bigr)\nonumber\\
&\leq& \lim_{n\rightarrow+\infty}
\PP\Bigl(\inf_{t\in[0,T\wedge \bar t_{\varepsilon'}]} \NISn_{t\wedge
\tau^n_\varepsilon}>\varepsilon\Bigr)\\
&=& \lim_{n\rightarrow+\infty}
\PP(\tau^n_\varepsilon> T\wedge
\bar t_{\varepsilon'}).\nonumber
\end{eqnarray}
%
Hence, we have
\[
\Psi_{\cdot\wedge\tau^n_\varepsilon\wedge \bar t_{\varepsilon'}\wedge
T}^{{\II}\sscS,f}\bigl(\mu^{(n)}\bigr) = \Psi_{\cdot\wedge\tau
^n_\varepsilon
\wedge
T}^{{\II}\sscS,f}\bigl(\mu^{(n)}\bigr)\ind_{\tau^n_\varepsilon\leq
\bar t_{\varepsilon
'}\wedge T} +
\Psi^{{\II}\sscS,f}_{\cdot\wedge \bar t_{\varepsilon'}\wedge T}\bigl(\mu
^{(n)}_{\cdot\wedge\tau^n_\varepsilon}\bigr) \ind_{\tau^n_\varepsilon
>\bar t_{\varepsilon'}\wedge T}.
\]
From the estimates of the different terms in (\ref{decompo}), $\Psi
_{\cdot\wedge\tau^n_\varepsilon\wedge T}^{{\II}\sscS,f}(\mu
^{(n)})$ is
upper bounded by a moment of $\mu^{(n)}$ of order 4. In view of (\ref
{eqhypmoments}) and (\ref{etapelimitetaun}),
the first term in the right-hand side converges in $L^1$ and hence, in
probability, to zero. Using the continuity of $\Psi^{{\II}\sscS, f}$ on
$\D(\R_+,\M_{0,A} \times\M_{\varepsilon,A} \times\M_{0,A})$,
$\Psi^{{\II}\sscS,f}(\mu^{(n)}_{\cdot\wedge\tau^n_\varepsilon
})$ converges
to $\Psi^{{\II}\sscS,f}(\bar{\mu})$ and, therefore, $\Psi^{{\II
}\sscS
,f}_{\cdot\wedge \bar t_{\varepsilon'}\wedge T}(\mu^{(n)}_{\cdot\wedge
\tau
^n_\varepsilon})$ converges to $\Psi^{{\II}\sscS
,f}_{\cdot\wedge \bar t_{\varepsilon'}\wedge T}(\bar{\mu})$. Thanks to this
and (\ref{etapelimitetaun}), the second term in the right-hand side
converges to $\Psi_{\cdot\wedge \bar t_{\varepsilon
'}\wedge
T}^{{\II}\sscS,f}(\bar{\mu})$ in $\D(\R_+,\R)$.

Then, 
$ (\langle\muISn_{\cdot\wedge\tau^n_\varepsilon\wedge
\bar t_{\varepsilon
'}\wedge T},f\rangle- \Psi^{{\II}\sscS,f}_{\cdot\wedge\tau
_\varepsilon
^n\wedge \bar t_{\varepsilon'}\wedge T}(\mu^{(n)}))_{n\in\N^*}$
converges in probability to $\langle\bar{\mu}_{\cdot\wedge
\bar t_{\varepsilon'}\wedge T},f\rangle- \Psi^{{\II}\sscS,f}_{\cdot
\wedge
\bar t_{\varepsilon'}\wedge T}(\bar{\mu})$. From (\ref{decompo}), this
sequence also converges in probability to zero.
By identification of these limits, $\bmuIS$ solves (\ref{limitereseauinfiniSI}) on
$[0,\bar{t}_{\varepsilon'}\wedge T]$. If $\langle
\bmuRS_0,\chi\rangle>0$ then similar techniques can be used. Else, the
result is obvious since for all $t\in [0,\bar t_{\varepsilon'}\wedge
T]$, $\langle \muISn_t,\chi\rangle>\varepsilon$ and the term
$p^n_t(j,\ell|k-1)$ is negligible when $\ell>0$. Thus $\bar{\mu}$
coincides a.s. with the only continuous deterministic solution of
(\ref{limitereseauinfiniS})--(\ref{limitereseauinfiniSR}) on $[0,\bar
t_{\varepsilon'}\wedge T]$. This implies that $\bar
t_{\varepsilon'}=t_{\varepsilon'}$ and yields the convergence in
probability of $(\mu^{(n)}_{\cdot\wedge \tau^n_\varepsilon})_{n\in \N^*}$
to $\bar{\mu}$, uniformly on $[0,t_{\varepsilon'}]$ since $\bar{\mu}$
is continuous.

We finally prove that the nonlocalized sequence $(\mu^{(n)})_{n\in
\N^*}$ also converges uniformly and in probability to $\bar{\mu}$ in
$\D([0,t_{\varepsilon'}],\mathcal{M}_{0,A}\times\M_{\varepsilon
,A} \times
\M_{0,A})$. For a~small positive $\eta$,
%
%
\begin{eqnarray}\label{etape6}
&&\PP\Bigl( \sup_{t \in[0,t_{\varepsilon'}]}\bigl|\bigl\langle\muISn_{t},f\bigr\rangle-
\Psi^{{\II}\sscS,f}_{t}(\bar{\mu})\bigr|>\eta\Bigr)\nonumber\\[-3pt]
&&\qquad\leq\PP\biggl(\sup_{t\in[0,t_{\varepsilon'}]}\bigl|
\Psi^{{\II}\sscS,f}_{t\wedge\tau^n_\varepsilon}\bigl(\mu^{(n)}\bigr) -
\Psi^{{\II}\sscS,f}_{t
}(\bar{\mu})\bigr|>\frac{\eta}{2} ;
\tau^n_\varepsilon\geq t_{\varepsilon'} \biggr)\\[-3pt]
&&\qquad\quad{}+
\PP\biggl(\sup_{t\in
[0,t_{\varepsilon'}]}\bigl|\Delta^{n,f}_{t\wedge
\tau^n_\varepsilon}+M_{t\wedge
\tau^n_\varepsilon}^{(n),{\II}\sscS,f}\bigr|>
\frac{\eta}{2}\biggr)+\PP(\tau^n_\varepsilon<
t_{\varepsilon'}).
\nonumber
\end{eqnarray}
Using the continuity of $\Psi^f$ and the uniform convergence in
probability~pro\-ved above, the first term in the right-hand side of
(\ref{etape6}) converges to zero. We can show that the second term
converges to zero by using Doob's inequality together with the
estimates of the bracket of $M^{(n),{\II}\sscS,f}$ [similar to (\ref
{etape3})] and of $\Delta^{n,f}$ (Step~2). Finally, the third term
vanishes in view of~(\ref{etapelimitetaun}).

The convergence of the original sequence $(\mu^{(n)})_{n\in
\N^*}$ is then entailed by the uniqueness of the solution to
(\ref{limitereseauinfiniS})--(\ref{limitereseauinfiniSR}), implied by
Step 2.

\textit{Step} 5. When $n\rightarrow+\infty$, by taking the
limit in (\ref{eqmuSrenorm}), $(\muSn)_{n\in\N^*}$
converges in $\D(\R_+,\mathcal{M}_{0,A})$ to the solution of the following
transport equation, that can be solved in function of $\bpI$. For
every bounded function $f \dvtx (k,t)\mapsto f_t(k)\in
\Co_b^{0,1}(\N\times\R_+,\R)$ of class $\Co^1$ with bounded
derivative with respect to $t$,
%
%
\begin{equation}
\langle\bar{\mu}^{\sscS}_t,f_t\rangle= \langle
\bar{\mu}_0^{\sscS},f_0\rangle-\int_0^t \langle\bar{\mu}^{\sscS
}_s,r\chi\bpI_s
f_s - \partial_sf_s \rangle\,\rmd {s}.\vadjust{\goodbreak}
\end{equation}
Choosing $f(k,s)=\varphi(k)\exp(-rk\int_0^{t-s}\bpI(u)\,\rmd u)$, we
obtain that
%
%
\begin{equation}
\langle\bar{\mu}^{\sscS}_t,\varphi\rangle= \sum_{k\in\N}
\varphi
(k)\theta_t^k \bar{\mu}_0^{\sscS}(k),
\end{equation}
where $\theta_t=\exp(-r\int_0^{t}\bpI(u)\,\rmd u)$ is the probability
that a given degree 1 node remains susceptible at time $t$. This is the
announced equation (\ref{limitereseauinfiniS}).
\end{pf}

We end this section with a lower bound of the time $t_{\varepsilon'}$
until which we proved that the convergence to Volz's equations holds.
\begin{prop}
\label{prohorizon}
Under the assumptions of Theorem \ref{propconvergencemunS}, 
%
%
\begin{equation}
\label{eqdeftauepsilon}
t_{\varepsilon'}>\bar\tau_{\varepsilon'}:=\frac{\log(\langle
\bmuS _0,\chi^2 \rangle+\bNIS_0)-\log(\langle\bmuS_0,\chi^2
\rangle+\varepsilon
')}{\max(\beta, r)}.
\end{equation}
\end{prop}
\begin{pf}
Because of the moment assumption (\ref{eqhypmoments}), we can
prove\break
that~(\ref{decompo}) also holds for $f=\chi$. This is obtained by
replacing in (\ref{etape777}), (\ref{etape5}), (\ref{etape753}) and
(\ref{etape779}) $\|f\|_\infty$ by $k$ and using the assumption of
boundedness of the moments of order 5 in (\ref{etape5}) and (\ref
{etape779}). This shows that $(\NISn)_{n\in\N}$ converges, uniformly
on $[0,t_{\varepsilon'}]$ and in
probability, to the deterministic and continuous solution
$\bNIS=\langle\bmuIS,\chi\rangle$. We introduce the event
$\mathcal A^{n}_{\xi}= \{\mid\NISnn_0-n \bNIS_0\mid\le\xi\}$
where their differences are bounded by $\xi>0$.
Recall the definition (\ref{eqdeftauepsilonn}) and let us introduce
the number of edges $Z^n_t$ that were $\cI\cS$ at time $0$ and that
have been removed before~$t$. For $t\geq\tau^n_{\varepsilon'}$, we
have necessarily that $Z^n_t\geq N^{n,\II\sscS}_0-n\varepsilon'$. Thus,
%
%
\begin{eqnarray}\label{etape12}
{\mathbb P}(\{\tau^n_{\varepsilon'} \le t\}\cap\mathcal A^n_\xi
) &\leq& {\mathbb P}(\{Z^n_t > \NISnn_0 - n\varepsilon' \} \cap
\mathcal A^n_{\xi})\nonumber\\[-8pt]\\[-8pt]
&\le&{\mathbb P}\bigl(\{Z^n_t > n(\bNIS_0 -\varepsilon')-\xi\} \cap
\mathcal A^n_{\xi}\bigr).\nonumber
\end{eqnarray}
When susceptible (resp., infectious) individuals of degree $k$ are
contaminated (resp., removed), at most
$k$ $\cI\cS$-edges are lost. Let $X^{n,k}_t$ be the number of edges
that, at time $0$, are $\cI\cS$ with susceptible alter of degree $k$,
and that have transmitted the disease before time $t$. Let $Y^{n,k}_t$
be the number of initially infectious individuals $x$ with $d_x(\cS
_0)=k$ and who have been removed before time $t$. $X^{n,k}_t$ and
$Y^{n,k}_t$ are bounded by $k \muSnn_0(k)$ and $\muSInn_0(k)$. Thus,
%
%
\begin{equation}\label{eqmajoreZ1}
Z^n_t\leq\sum_{k\in\N}k(X^{n,k}_t+ Y^{n,k}_t).
\end{equation}
%
Let us stochastically upper bound $Z^n_t$. 
Since each $\cI\cS$-edge transmits the disease independently at rate
$r$, $X^{n,k}_t$ is stochastically dominated by a binomial r.v. of parameters
$k \muSnn_0(k)$ and $1-e^{-rt}$. We proceed similarly for~$Y^{n,k}_t$.
Conditionally to the initial condition, $X^{n,k}_t+Y^{n,k}_t$ is thus
stochastically dominated by a binomial r.v. $\tilde Z_t^{n,k}$ of
parameters $(k \muSnn_0(k)+ \muISnn_0(k))$ and $1-e^{-\max(\beta
,r)t}$. Then (\ref{etape12}) and (\ref{eqmajoreZ1}) give
%
%
\begin{equation}\label{eqhorizon1}
{\mathbb P}(\{\tau^n_{\varepsilon'} \le t\}\cap\mathcal A^n_\xi
)\leq{\mathbb P}\biggl(\sum_{k\in\N} \frac{k \tilde Z_t^{n,k}}{n}> \bNIS
_0 -\varepsilon'-\frac{\xi}{n}\biggr).
\end{equation}
Thanks\vspace*{1pt} to Assumption \ref{hypconvcondinit} and (\ref{eqhypmoments}),
the series $\sum_{k\in\N} k \tilde Z_t^{n,k}/n$ converges in $L^1$
and hence, in probability to $(\langle\bmuS_0,\chi^2\rangle+\bNIS
_0)(1-e^{-\max(\beta,r)t})$ when $n\rightarrow+\infty$. Thus, for
sufficiently large $n$,
\[
{\mathbb P}(\{\tau^n_{\varepsilon'} \le t\}\cap\mathcal A^n_\xi)=
1 \qquad\mbox{if }t>\bar{\tau}_{\varepsilon'}\quad\mbox{and}\quad 0\qquad \mbox{if
}t<\bar{\tau}_{\varepsilon'}.
\]
For all $t<\bar{\tau}_{\varepsilon'}$, it follows from Assumption
\ref{hypconvcondinit}, (\ref{eqhypmoments}) and Lemma \ref
{lemmeconvfnonbornee} that
\[
\lim_{n\rightarrow+\infty}{\mathbb P}(\tau^n_{\varepsilon'} \le t)
\le
\lim_{n\rightarrow+\infty}\bigl({\mathbb P}(\{\tau^n_{\varepsilon'} \le
t\} \cap\mathcal A^n_{\xi})+\PP((\mathcal A^n_\xi)^c)\bigr)
= 0,
\]
so that by Theorem \ref{propconvergencemunS}
\[
1= \lim_{n\rightarrow+\infty}{\mathbb P}(\tau^n_{\varepsilon'}
\geq\bar{\tau}_{\varepsilon'}) = \lim_{n\rightarrow
+\infty}\PP\Bigl(\inf_{t\leq\bar{\tau}_{\varepsilon'}}\NISn_t \geq
\varepsilon'\Bigr) = \PP\Bigl(\inf_{t\leq
\bar{\tau}_{\varepsilon'}}\bNIS_t \geq\varepsilon'\Bigr).
\]
This shows that $t_{\varepsilon'} \ge\bar\tau_{\varepsilon'}$
a.s., which concludes the proof.
\end{pf}

\subsection{Proof of Volz's equations}

\begin{prop}\label{propvolz}
The system (\ref{limitereseauinfiniS})--(\ref{limitereseauinfiniSR})
implies Volz's equations (\ref{volz1})--(\ref{volz4}).
\end{prop}

Before proving
Proposition \ref{propvolz}, we begin with a corollary of Theorem \ref
{propconvergencemunS}.
\begin{corol}\label{corolnbrearete}
For all $t\in \R_+$
%
%
\begin{eqnarray}\label{eqcorol}
\bar{N}_t^{\sscS} &=& \theta_t g'(\theta_t),\nonumber\\
\bar{N}_t^{{\II}\sscS} &=& \bar{N}_0^{{\II}\sscS}+\int_0^t r \bpI
_s \theta
_s g'(\theta_s)\biggl((\bpS_s-\bpI_s)\theta_s\frac{g''(\theta
_s)}{g'(\theta_s)}-1\biggr)-\beta\bar{N}_s^{{\II}\sscS} \,\rmd {s},\hspace*{-30pt}\\
\bar{N}_t^{\sRr\sscS}&=&\int_0^t \bigl(\beta
\bar{N}_s^{{\II}\sscS}-r\bpR_s\bpI_s\theta_s^2
g''(\theta_s)\bigr)\,\rmd {s}.\nonumber
\end{eqnarray}
\end{corol}
\begin{pf}
In the proof of Proposition \ref{prohorizon}, we have shown that
$(\NISn)_{n\in\N}$ converges uniformly on compact intervals
and in
probability to the deterministic and continuous solution
$\bNIS=\langle\bmuIS,\chi\rangle$.
Equation (\ref{limitereseauinfiniS}) with $f=\chi$ reads
%
%
\begin{equation}\label{etape9}
\bar{N}_t^{\sscS}=\sum_{k\in\N}\bar{\mu}_0^{\sscS}(k) k \theta
_t^k=\theta_t
\sum_{k=1}^{+\infty}\bar{\mu}_0^{\sscS}(k) k\theta_t^{k-1}=\theta_t
g'(\theta_t),
\end{equation}
that is, the first assertion of (\ref{eqcorol}).\vadjust{\goodbreak}

Choosing $f=\chi$ in (\ref{limitereseauinfiniSI}), we obtain
\begin{eqnarray*}
\bNIS_t &=& \bNIS_0-\int_0^t \beta\bNIS_s\,\rmd {s}
+ \int_0^t \sum_{k\in\N} \lambda_s(k) \sum_{j+\ell\leq
k-1}(k-2j-2-\ell)\\
&&{}\times\biggl[\frac{(k-1)!}{j!(k-1-j-\ell)!\ell!}(\bpI_s)^{j}(\bpR
_s)^{\ell}(\bpS_s)^{k-1-j-\ell}
\biggr]
\bmuS_s(k)\,\rmd {s}.
\end{eqnarray*}
Notice that the term in the square brackets is the probability to
obtain $(j,\ell,k-1-j-\ell)$ from a draw in
the multinomial distribution of parameters
$(k-1$, $(\bpI_s,\bpR_s,\bpS_s))$. Hence,
\[
\sum_{j+\ell\leq k-1}j\times
\biggl(\frac{(k-1)!}{j!(k-1-j-\ell)!\ell!}(\bpI_s)^{j}(\bpR_s)^{\ell
}(\bpS_s)^{k-1-j-\ell}
\biggr) =(k-1)\bpI_s
\]
as we recognize the
mean number of edges to $\I_s$ of an individual of degree
$k$. Other terms are treated similarly. Hence, with the definition of
$\lambda_s(k)$,
(\ref{eqdeflambda}),
\begin{eqnarray*}
\bNIS_t &=& \bNIS_0+ \int_0^t
r \bar{p}_s^{{\II}}\bigl(\langle\bmuS_s,\chi^2-2\chi\rangle
-(2\bpI_s+\bpR_s) \langle
\bmuS_s,\chi(\chi-1)\rangle\bigr) \,\rmd {s}\\
&&{} - \int_0^t \beta
\bNIS_s\,\rmd {s}.
\end{eqnarray*}
But since
\begin{eqnarray*}
\langle\bar{\mu}_t^{\sscS},\chi(\chi-1)\rangle &=& \sum_{k\in\N
}\bar{\mu}_0^{\sscS}(k) k(k-1)\theta_t^k = \theta_t^2 g''(\theta
_t),\\
\langle\bar{\mu}_t^{\sscS},\chi^2-2\chi\rangle &=& \langle
\bar{\mu}_t^{\sscS},\chi(\chi-1)\rangle- \langle\bar{\mu}_t^S,\chi\rangle=
\theta_t^2 g''(\theta_t)-\theta_t g'(\theta_t),
\end{eqnarray*}
we obtain by noticing that
$1-2\bpI_s-\bpR_s=\bpS_s-\bpI_s$,
%
%
\begin{equation}
\bNIS_t= \bNIS_0+\int_0^t
r \bpI_s\bigl( (\bpS_s-\bpI_s)
\theta_s^2g''(\theta_s)-\theta_s g'(\theta_s) \bigr) \,\rmd {s}-\int_0^t
\beta\bNIS_s\,\rmd {s},\hspace*{-35pt}
\end{equation}
which is the second assertion of (\ref{eqcorol}). The
third equation of (\ref{eqcorol}) is obtained similarly.
\end{pf}

We are now ready to prove Volz's equations.
\begin{pf*}{Proof of Proposition \ref{propvolz}}
We begin with the proof of (\ref{volz1}) and~(\ref{volz2}). Fix again
$t\ge0$. For the
size of the susceptible population, taking $\varphi=\car$ in~(\ref{limitereseauinfiniS}),
we are led to introduce the same
quantity $\theta_t=\exp(-r\int_0^t \bpI_s \,\rmd s)$ as Volz and obtain~(\ref{volz1}). For the size
of the infective population, setting $f=\car$ in (\ref
{limitereseauinfiniSI}) entails
\begin{eqnarray*}
\bar{I}_t &=& \bar{I}_0+\int_0^t \biggl(\sum_{k\in\N} rk\bpI_s \bar
{\mu}_s^{\sscS}(k)-\beta\bar{I}_s\biggr) \,\rmd {s}\\
&=& \bar{I}_0+\int_0^t \biggl( r\bpI_s \sum_{k\in\N}
\bar{\mu}_0^{\sscS}(k) k \theta_s^k-\beta\bar{I}_s\biggr) \,\rmd {s}\\
&=&
\bar{I}_0+\int_0^t \bigl( r\bpI_s \theta_s g'(\theta_s)-\beta
\bar{I}_s\bigr) \,\rmd {s}
\end{eqnarray*}
by using (\ref{limitereseauinfiniS}) with $f=\chi$ for the second
equality.

Let us\vspace*{1pt} now consider the probability that an edge with a susceptible
ego has an infectious alter. Both equations
(\ref{volz1}) and (\ref{volz2}) depend on $\bpI_t=\bNIS_t/\bNS_t$.
It is thus
important to obtain an equation for this quantity. In \cite{volz},
this equation also leads to introduce the quantity
$\bpS_t$.

From Corollary \ref{corolnbrearete}, we see that $\bNS$ and $\bNIS$
are differentiable and
\begin{eqnarray*}
\frac{\rmd \bpI_t}{\rmd {t}} &=& \frac{
\rmd}{\rmd {t}}\biggl(\frac{\bNIS_t}{\bNS_t}\biggr)=\frac{1}{\bNS_t}\,\frac
{\rmd }{\rmd {t}}(\bNIS_t)-\frac{\bNIS_t}{(\bNS
_t)^2}\,\frac{\rmd }{\rmd {t}}(\bNS_t)\\
&=& \biggl(r\bpI_t(\bpS_t-\bpI_t)\theta_t \frac{g''(\theta
_t)}{g'(\theta_t)}-r\bpI_t-\beta\bpI_t\biggr)\\
&&{} -\biggl(\frac{\bpI_t}{\theta_t g'(\theta_t)}\bigl(-r\bpI_t
\theta_t g'(\theta_t)+\theta_t g''(\theta_t)(-r\bpI_t \theta
_t)\bigr)\biggr)\\
&=& r\bpI_t\bpS_t \theta_t \frac{g''(\theta_t)}{g'(\theta_t)}
-r\bpI_t(1-\bpI_t)-\beta\bpI_t
\end{eqnarray*}
by using the equations 1 and 2 of (\ref{eqcorol}) for the derivatives
of $\bNS$ and $\bNIS$ with respect to time for the second line. This
achieves the proof of (\ref{volz3}).

For (\ref{volz4}), we notice
that $
\bpS_t=1-\bpI_t-\bpR_t$ and achieve the proof by showing that
%
%
\begin{equation}
\bpR_t= \int_0^t (\beta\bpI_s-r\bpI_s
\bpR_s )\,\rmd {s}
\end{equation}
by using arguments similar as for $\bpI_t$.
\end{pf*}
\begin{Rque}Miller \cite{miller} shows that Volz's equations can be
reduced to only three ODEs:
\begin{eqnarray*}
\bar{S}_t&=&g(\theta_t),\qquad \frac{\rmd \bar{R}_t}{
\rmd{t}}=\beta\bar{I}_t,\qquad \bar{I}_t=(\bar{S}_0+\bar{I}_0)-\bar
{S}_t-\bar{R}_t,\\
\frac{\rmd \theta_t}{\rmd {t}}&=&-r\theta_t+\beta
(1-\theta_t)+\beta\frac{g'(\theta_t)}{g'(1)}.
\end{eqnarray*}
The last ODE is obtained by considering the probability that an edge
with an infectious ego drawn at random has not transmitted the disease.
However, in his simplifications, he uses that the degree distributions
$\bmuS_0/\bar{S}_0$ and $\sum_{k\in\N}p_k \delta_k$ are the same,
which is not necessarily the case\vspace*{1pt} (see our
Remark~\ref{rqueconvcondinit}). Moreover, it is more natural to have an ODE
on $\bar{I}_t$ and $\bNIS_t$ is a~natural quantity that is of
interest in itself for the dynamics.\vadjust{\goodbreak}
\end{Rque}

\begin{appendix}\label{appendiceA}

\section*{Appendix: Finite measures on $\N$}

First, some notation is needed in order to clarify the way the atoms of
a given element of $\M_F(\N)$ are ranked.
For all $\mu\in\M_F(\N)$, let $F_\mu$ be its cumulative
distribution function and $F_\mu^{-1}$ be its right inverse defined as
%
%
\setcounter{equation}{0}
\begin{equation}\label{etape7}
\forall x\in\R_+\qquad F^{-1}_{\mu}(x)=\inf\{i\in\N, F_{\mu}(i)\geq
x\}.
\end{equation}
Let $\mu=\sum_{n\in\N}a_n\delta_{n}$ be an integer-valued measure
of $\M_F(\N)$, that is, such that the $a_n$'s are integers
themselves. Then, for each atom
$n\in\N$ of $\mu$ such that $a_n>0$, we duplicate the atom $n$ with
multiplicity $a_n$, and we rank the atoms of $\mu$ by increasing
values, sorting arbitrarily the atoms having the same value. Then, we
denote for any $i \le\langle\mu,\mathbf1 \rangle$,
%
%
\begin{equation}
\label{eqnumeromu}
\gamma_i(\mu)=F_{\mu}^{-1}(i),
\end{equation}
the level of the $i$th atom of the measure, when ranked as described
above. We refer to Example \ref{exempleFmu} for a simple illustration.

We now make precise a few topological properties of spaces of measures
and measure-valued processes.
For $T>0$ and a Polish space $(E,d_E)$, we denote by $\D([0,T],E)$ the
Skorokhod space of c\`{a}dl\`{a}g (right-continuous left-limited)
functions from $\R$ to $E$ (e.g., \cite{billingsley,joffemetivier})
equipped with the Skorokhod topology induced by the metric
%
%
\begin{eqnarray}\label{distSkorokhod}
d_T(f,g)&:=&\inf_{\alpha\in\Delta([0,T])} \biggl\{\mathop{\sup_{(s,t)\in
[0,T]^2,}}_{s \ne t} \biggl| \log\frac{\alpha(s)-\alpha(t)}{s-t}\biggr| \nonumber\\[-8pt]\\[-8pt]
&&\hspace*{59.7pt}{} + \sup
_{t\le T} d_E( f(t),g(\alpha(t)) ) \biggr\},\nonumber
\end{eqnarray}
where the infimum is taken over the set $\Delta([0,T])$ of continuous
increasing functions $\alpha\dvtx [0,T] \to[0,T]$ such that $\alpha
(0)=0$ and $\alpha(T)=T$.

Limit theorems are heavily dependent on the topologies considered.
We introduce here several technical lemmas on the space of measures
related to these questions. For any fixed $0 \le\varepsilon< A$,
recall the definition of $\mathcal{M}_{\varepsilon,A}$ in~(\ref{eqdefM}).
Remark that for any $\nu\in\M_{\varepsilon, A}$, and $i\in\{
0,\ldots,5\}$,
$\langle\nu,\chi^i \rangle\leq A$ since the support of $\nu$ is
included in
$\N$.
\setcounter{prop}{0}
\begin{lemme}
\label{lemuniformintegrability}
Let ${\mathfrak I}$ a set and a family $(\nu_\tau,\tau\in{\mathfrak
I})$ of elements of $\M_{\varepsilon, A}$. Then, for any real
function $f$ on $\N$ such that
$f(k)=o(k^5)$, we have that
\[
\lim_{K\to\infty} \sup_{\tau\in{\mathfrak I}}\bigl|\bigl\langle\nu_\tau,
f\car_{[K,\infty)}\bigr\rangle\bigr|=0.
\]
%
\end{lemme}
\begin{pf}
By the Markov inequality, for any $\tau\in{\mathfrak I}$, for any $K$,
we have
\[
\sum_{k\ge K} |f(k)|\nu_\tau(k)\le A \sup_{k\ge K}
\frac{|f(k)|}{k^5},\vadjust{\goodbreak}
\]
hence,
\[
\lim_{K\to\infty} \sup_{\tau\in{\mathfrak I}}|\langle\nu_\tau,
f\rangle|\le A \limsup_{k\to\infty} \frac{|f(k)|}{k^5}=0.
\]
The proof is thus complete.
\end{pf}
\begin{lemme}
\label{lemmapropM}
For any $A>0$, the set $\M_{\varepsilon, A}$ is a closed subset of
$\M_F(\N)$
embedded with the topology of weak convergence.
\end{lemme}
\begin{pf}
Let $(\mu_n)_{n\in\N}$ be a sequence of $\M_{\varepsilon, A}$
converging to $\mu
\in\M_F(\N)$ for the weak topology, which implies in particular
that $\lim_{n\rightarrow+\infty}\mu_n(k)= \mu(k)$ for any
$k\in\N$. Denoting for all $n$ and $k\in\N$, $f_n(k)=k^5\mu^n(k)$,
we have that $\lim_{n\rightarrow+\infty} f_n(k)= f(k):=k^5\mu(k)$,
$\mu$-a.e., and Fatou's lemma implies
\[
\langle\mu,\chi^5 \rangle=\sum_{k\in\N} f(k)\le\lim\inf
_{n\rightarrow
\infty}\sum_{k\in\N} f_n(k)=\lim\inf_{n\rightarrow\infty
}\langle\mu_n,\chi^5 \rangle.
\]
Since $\langle\mu_n,\mathbf1 \rangle$ tends to $\langle\mu
,\mathbf1 \rangle$, we
have that $\langle\mu,\car+ \chi^5 \rangle \le A$.

Furthermore, by uniform integrability (Lemma \ref
{lemuniformintegrability}), it is also clear that
\[
\varepsilon\le\lim_{n\to\infty} \langle\mu_n, \chi \rangle
=\langle\mu, \chi \rangle,
\]
which shows that $\mu\in\M_{\varepsilon, A}$.
\end{pf}
\begin{lemme}
\label{lemmatotalvariation}
The traces on $\M_{\varepsilon, A}$ of the total variation topology
and of the weak
topology coincide.
\end{lemme}
\begin{pf}
It is well known that the total variation topology is coarser than
the weak topology. In the reverse direction, assume that
$(\mu_n)_{n\in\N}$ is a~sequence of weakly converging measures
belonging to $\M_{\varepsilon, A}$. Since
\[
d_{\mathrm{TV}}(\mu_n, \mu)\le\sum_{k\in\N} |\mu_n(k)-\mu(k)|
\]
according to Lemma \ref{lemuniformintegrability}, it is then easily
deduced that the right-hand side converges to $0$ as $n$ goes to infinity.
\end{pf}
\begin{lemme}\label{lemmeconvfnonbornee}If the sequence
$(\mu_n)_{n\in\N}$ of $\M_{\varepsilon, A}^\N$ converges weakly
to the measure
$\mu\in\M_{\varepsilon, A}$, then $(\langle\mu_n,f\rangle)_{n\in
\N}$ converges
to $\langle\mu,f\rangle$ for all function $f$ such that
$f(k)=o(k^5)$ for all large $k$.
\end{lemme}
\begin{pf} Triangular inequality says that
\begin{eqnarray*}
|\langle\mu_n,f\rangle- \langle\mu,f\rangle|&\leq&\bigl|\bigl\langle\mu
_n,f \ind_{[0,K]}\bigr\rangle- \bigl\langle\mu,f\ind_{[0,K]}\bigr\rangle\bigr|\\
&&{} +\bigl|\bigl\langle\mu,f\ind_{(K,+\infty)}\bigr\rangle\bigr|+\bigl|\bigl\langle\mu_n,f\ind
_{(K,+\infty)}\bigr\rangle\bigr|.
\end{eqnarray*}
We then conclude by uniform integrability and weak convergence.\vadjust{\goodbreak}
\end{pf}
Recall that $\M_{\varepsilon, A}$ can be embedded with the total
variation distance
topology, hence, the topology on $\D([0,T],\M_{\varepsilon, A})$ is
induced by the
distance
\[
\rho_T(\mu_{\cdot}, \nu_{\cdot})=\inf_{\alpha\in
\Delta([0,T])}\biggl(\mathop{\sup_{(s,t)\in[0,T]^2,}}_{s \ne
t} \biggl| \log\frac{\alpha(s)-\alpha(t)}{s-t}\biggr|
+\sup_{t\le T} d_{\mathrm{TV}}\bigl(\mu_t, \nu_{\alpha(t)}\bigr)\biggr).
\]
\begin{lemme}
\label{lemmaconti1}
For any $p\le5$, the following map is continuous:
\[
\Phi_p \dvtx \cases{
\D(\R_+,\M_{\varepsilon, A})\longrightarrow
\D(\R_+,\R),\cr
\nu_{\cdot}\longmapsto\langle\nu_{\cdot}, \chi^p
\rangle.}
\]
%
\end{lemme}
\begin{pf}It is sufficient to prove the continuity of the above
mappings from $\D([0,T],\M_{\varepsilon, A})$ to $\D([0,T],\R)$,
for any $T \ge0$, where the latter are equipped with the corresponding
Skorokhod topologies.
For $\mu$ and $\nu$ two elements of $\M_{\varepsilon, A}$, for any
$p\le5$, for any
positive integer $K$,
according to the Markov inequality,
%
%
\begin{eqnarray}\label{eq1}
|\langle\mu, \chi^p\rangle- \langle\nu, \chi^p\rangle| &\le&
2\frac A{K^p} + \bigl|\bigl\langle\mu-\nu, \chi^p\ind_{[0, K]}\bigr\rangle\bigr|\nonumber\\[-8pt]\\[-8pt]
&\le&
2\frac A{K^p} +K^p d_{\mathrm{TV}}(\mu, \nu).\nonumber
\end{eqnarray}
Using (\ref{distSkorokhod}) and (\ref{eq1}) we have for any $K>0$,
%
\[
d_T(\langle\mu_{\cdot}, \chi^p \rangle, \langle\nu_{\cdot},
\chi^p \rangle)\le2\frac
A{K^p} + K^p d_T(\mu_{\cdot}, \nu_{\cdot}),
\]
and hence, the continuity of $\Phi_p$.
\end{pf}
\end{appendix}

\section*{Acknowledgments}

Tran thanks T. L. Parsons for the invitation to the \mbox{DIMACS}
\textit{Workshop on Stochasticity in Population and Disease Dynamics}
in December 2008 and L. M. Wahl for discussions at this workshop on
epidemics on graphs.
The authors also thank M. Costa and E. Pardoux for their careful
reading and for discussions which improved the manuscript.


%

%
\printaddresses


\begin{thebibliography}{29}

\bibitem{andersson}
%
\begin{barticle}[author]
\bauthor{\bsnm{Andersson},~\bfnm{H.}\binits{H.}}
(\byear{1998}).
\btitle{Limit theorems for a random graph epidemic model}.
\bjournal{Ann. Appl. Probab.}
\bvolume{8}
\bpages{1331--1349}.
\end{barticle}
%
\MR{1661200}
\endbibitem

\bibitem{andersonbritton}
%
\begin{bbook}[author]
\bauthor{\bsnm{Andersson},~\bfnm{H.}\binits{H.}} \AND
\bauthor{\bsnm{Britton},~\bfnm{T.}\binits{T.}}
(\byear{2000}).
\btitle{Stochastic Epidemic Models and Their Statistical Analysis}.
\bseries{Lecture Notes in Statistics}
\bvolume{151}.
\bpublisher{Springer}, \baddress{New York}.
\end{bbook}
%
\MR{1784822}
\endbibitem

\bibitem{ballneal}
%
\begin{barticle}[author]
\bauthor{\bsnm{Ball},~\bfnm{F.}\binits{F.}} \AND
\bauthor{\bsnm{Neal},~\bfnm{P.}\binits{P.}}
(\byear{2008}).
\btitle{Network epidemic models with two levels of mixing}.
\bjournal{Math. Biosci.}
\bvolume{212}
\bpages{69--87}.
\end{barticle}
%
\MR{2399833}
\endbibitem

\bibitem{barthelemybarratpastorsatorrasvespignani}
%
\begin{barticle}[author]
\bauthor{\bsnm{Barth{\'{e}}lemy},~\bfnm{M.}\binits{M.}},
\bauthor{\bsnm{Barrat},~\bfnm{A.}\binits{A.}},
\bauthor{\bsnm{Pastor-Satorras},~\bfnm{R.}\binits{R.}} \AND
\bauthor{\bsnm{Vespignani},~\bfnm{A.}\binits{A.}}
(\byear{2005}).
\btitle{Dynamical patterns of epidemic outbreaks in complex heterogeneous
networks}.
\bjournal{J. Theoret. Biol.}
\bvolume{235}
\bpages{275--288}.
\end{barticle}
%
\MR{2157753}
\endbibitem

\bibitem{bartlett}
%
\begin{bbook}[author]
\bauthor{\bsnm{Bartlett},~\bfnm{M.~S.}\binits{M.~S.}}
(\byear{1960}).
\btitle{Stochastic Population Models in Ecology and Epidemiology}.
\bpublisher{Methuen}, \baddress{London}.
\end{bbook}
%
\MR{0118550}
\endbibitem

\bibitem{billingsley}
%
\begin{bbook}[author]
\bauthor{\bsnm{Billingsley},~\bfnm{P.}\binits{P.}}
(\byear{1968}).
\btitle{Convergence of Probability Measures}.
\bpublisher{Wiley}, \baddress{New York}.
\end{bbook}
%
\MR{0233396}
\endbibitem

\bibitem{bollobas2001}
%
\begin{bbook}[author]
\bauthor{\bsnm{Bollob{\'{a}}s},~\bfnm{B.}\binits{B.}}
(\byear{2001}).
\btitle{Random Graphs}, \bedition{2nd} ed.
\bpublisher{Cambridge Univ. Press}, \baddress{Cambridge}.
\end{bbook}
%
\MR{1864966}
\endbibitem

\bibitem{clemencondearazozarossitran}
%
\begin{bmisc}[author]
\bauthor{\bsnm{Cl{\'{e}}men{\c{c}}on},~\bfnm{S.}\binits{S.}},
\bauthor{\bparticle{De}~\bsnm{Arazoza},~\bfnm{H.}\binits{H.}},
\bauthor{\bsnm{Rossi},~\bfnm{F.}\binits{F.}} \AND
\bauthor{\bsnm{Tran},~\bfnm{V.~C.}\binits{V.~C.}}
\bhowpublished{A network analysis of the {H}{I}{V}--{A}{I}{D}{S}
epidemic in {C}uba. Unpublished manuscript.}
\end{bmisc}
%
\endbibitem

\bibitem{arazozaclemencontran}
%
\begin{barticle}[author]
\bauthor{\bsnm{Cl{\'{e}}men{\c{c}}on},~\bfnm{S.}\binits{S.}},
\bauthor{\bsnm{Tran},~\bfnm{V.~C.}\binits{V.~C.}} \AND
\bauthor{\bparticle{De}~\bsnm{Arazoza},~\bfnm{H.}\binits{H.}}
(\byear{2008}).
\btitle{A stochastic {S}{I}{R} model with contact-tracing: Large population
limits and statistical inference}.
\bjournal{J. Biol. Dyn.}
\bvolume{2}
\bpages{392--414}.
\end{barticle}
%
\endbibitem

\bibitem{durrett}
%
\begin{bbook}[author]
\bauthor{\bsnm{Durrett},~\bfnm{R.}\binits{R.}}
(\byear{2007}).
\btitle{Random Graph Dynamics}.
\bpublisher{Cambridge Univ. Press}, \baddress{Cambridge}.
\end{bbook}
%
\MR{2271734}
\endbibitem

\bibitem{ethierkurtz}
%
\begin{bbook}[author]
\bauthor{\bsnm{Ethier},~\bfnm{S.~N.}\binits{S.~N.}} \AND
\bauthor{\bsnm{Kurtz},~\bfnm{T.~G.}\binits{T.~G.}}
(\byear{1986}).
\btitle{Markov Processus, Characterization and Convergence}.
\bpublisher{Wiley}, \baddress{New York}.
\end{bbook}
%
\MR{0838085}
\endbibitem

\bibitem{evans}
%
\begin{bbook}[author]
\bauthor{\bsnm{Evans},~\bfnm{L.~C.}\binits{L.~C.}}
(\byear{1998}).
\btitle{Partial Differential Equations}.
\bseries{Graduate Studies in Mathematics}
\bvolume{19}.
\bpublisher{Amer. Math. Soc.}, \baddress{Providence, RI}.
\end{bbook}
%
\MR{1625845}
\endbibitem

\bibitem{fourniermeleard}
%
\begin{barticle}[author]
\bauthor{\bsnm{Fournier},~\bfnm{N.}\binits{N.}} \AND
\bauthor{\bsnm{M{\'{e}}l{\'{e}}ard},~\bfnm{S.}\binits{S.}}
(\byear{2004}).
\btitle{A microscopic probabilistic description of a locally regulated
population and macroscopic approximations}.
\bjournal{Ann. Appl. Probab.}
\bvolume{14}
\bpages{1880--1919}.
\end{barticle}
%
\MR{2099656}
\endbibitem

\bibitem{ikedawatanabe}
%
\begin{bbook}[author]
\bauthor{\bsnm{Ikeda},~\bfnm{N.}\binits{N.}} \AND
\bauthor{\bsnm{Watanabe},~\bfnm{S.}\binits{S.}}
(\byear{1989}).
\btitle{Stochastic Differential Equations and Diffusion Processes},
\bedition{2nd} ed.
\bseries{North-Holland Mathematical Library}
\bvolume{24}.
\bpublisher{North-Holland}, \baddress{Amsterdam}.
\end{bbook}
%
\MR{1011252}
\endbibitem

\bibitem{jacod}
%
\begin{bbook}[author]
\bauthor{\bsnm{Jacod},~\bfnm{J.}\binits{J.}} \AND
\bauthor{\bsnm{Shiryaev},~\bfnm{A.~N.}\binits{A.~N.}}
(\byear{1987}).
\btitle{Limit Theorems for Stochastic Processes}.
\bpublisher{Springer}, \baddress{Berlin}.
\end{bbook}
%
\MR{0959133}
\endbibitem

\bibitem{joffemetivier}
%
\begin{barticle}[author]
\bauthor{\bsnm{Joffe},~\bfnm{A.}\binits{A.}} \AND
\bauthor{\bsnm{M{\'{e}}tivier},~\bfnm{M.}\binits{M.}}
(\byear{1986}).
\btitle{Weak convergence of sequences of semimartingales with
applications to
multitype branching processes}.
\bjournal{Adv. in Appl. Probab.}
\bvolume{18}
\bpages{20--65}.
\end{barticle}
%
\MR{0827331}
\endbibitem

\bibitem{kermackmckendrick}
%
\begin{barticle}[author]
\bauthor{\bsnm{Kermack},~\bfnm{W.~O.}\binits{W.~O.}} \AND
\bauthor{\bsnm{McKendrick},~\bfnm{A.~G.}\binits{A.~G.}}
(\byear{1927}).
\btitle{A contribution to the mathematical theory of epidemics}.
\bjournal{Proc. R. Soc. Lond. Ser. A Math. Phys. Eng. Sci.}
\bvolume{115}
\bpages{700--721}.
\end{barticle}
%
\endbibitem

\bibitem{miller}
%
\begin{barticle}[author]
\bauthor{\bsnm{Miller},~\bfnm{J.~C.}\binits{J.~C.}}
(\byear{2011}).
\btitle{A note on a paper by {E}rik {V}olz: {S}{I}{R}
dynamics in random networks}.
\bjournal{J.~Math. Biol.}
\bvolume{62}
\bpages{349--358}.
\end{barticle}
%
\endbibitem

\bibitem{molloyreed}
%
\begin{barticle}[author]
\bauthor{\bsnm{Molloy},~\bfnm{M.}\binits{M.}} \AND
\bauthor{\bsnm{Reed},~\bfnm{B.}\binits{B.}}
(\byear{1995}).
\btitle{A critical point for random graphs with a given degree sequence}.
\bjournal{Random Structures Algorithms}
\bvolume{6}
\bpages{161--180}.
\end{barticle}
%
\MR{1370952}
\endbibitem

\bibitem{newman}
%
\begin{barticle}[author]
\bauthor{\bsnm{Newman},~\bfnm{M.~E.~J.}\binits{M.~E.~J.}}
(\byear{2002}).
\btitle{The spread of epidemic disease on networks}.
\bjournal{Phys. Rev. E (3)}
\bvolume{66}
\bpages{016128, 11}.
\end{barticle}
%
\MR{1919737}
\endbibitem

\bibitem{newmanSIAM}
%
\begin{barticle}[author]
\bauthor{\bsnm{Newman},~\bfnm{M.~E.~J.}\binits{M.~E.~J.}}
(\byear{2003}).
\btitle{The structure and function of complex networks}.
\bjournal{SIAM Rev.}
\bvolume{45}
\bpages{167--256}.
\end{barticle}
%
\MR{2010377}
\endbibitem

\bibitem{newmanstrogatzwatts}
%
\begin{barticle}[author]
\bauthor{\bsnm{Newman},~\bfnm{M.~E.~J.}\binits{M.~E.~J.}},
\bauthor{\bsnm{Strogatz},~\bfnm{S.~H.}\binits{S.~H.}} \AND
\bauthor{\bsnm{Watts},~\bfnm{D.~J.}\binits{D.~J.}}
(\byear{2001}).
\btitle{Random graphs with arbitrary degree distributions and their
applications}.
\bjournal{Phys. Rev. E (3)}
\bvolume{64}.
\end{barticle}
%
\endbibitem

\bibitem{pastorsatorrasvespignani}
%
\begin{binproceedings}[author]
\bauthor{\bsnm{Pastor-Satorras},~\bfnm{R.}\binits{R.}} \AND
\bauthor{\bsnm{Vespignani},~\bfnm{A.}\binits{A.}}
(\byear{2002}).
\btitle{Epidemics and immunization in scale-free networks}.
In \bbooktitle{Handbook of Graphs and Networks: From the Genome to the
Internet}
\bpages{113--132}.
\bpublisher{Wiley-VCH}, \baddress{Berlin}.
\end{binproceedings}
\endbibitem

\bibitem{roelly}
%
\begin{barticle}[author]
\bauthor{\bsnm{Roelly-Coppoletta},~\bfnm{S.}\binits{S.}}
(\byear{1986}).
\btitle{A criterion of convergence of measure-valued processes:
Application to
measure branching processes}.
\bjournal{Stochastics}
\bvolume{17}
\bpages{43--65}.
\end{barticle}
%
\MR{0878553}
\endbibitem

\bibitem{chithese}
%
\begin{bmisc}[author]
\bauthor{\bsnm{Tran},~\bfnm{V.~C.}\binits{V.~C.}}
(\byear{2007}).
\bhowpublished{Mod\`{e}les particulaires stochastiques pour des probl\`{e}mes
d'\'{e}volution adaptative et pour l'approximation de solutions statistiques.
Ph.D. thesis, Univ. Paris X---Nanterre. Available at
\texttt{\href{http://tel.archives-ouvertes.fr/tel-00125100}{http://tel.archives-}
\href{http://tel.archives-ouvertes.fr/tel-00125100}{ouvertes.fr/tel-00125100}}.}
\end{bmisc}
%
\endbibitem

\bibitem{vanderhofstad}
%
\begin{bmisc}[author]
\bauthor{\bparticle{van~der} \bsnm{Hofstad},~\bfnm{R.}\binits{R.}}
(\byear{2011}).
\bhowpublished{{Random graphs and complex networks}.
Lecture Notes. To
appear. 
Available at \url{http://www.win.tue.nl/\textasciitilde rhofstad}.}
\end{bmisc}
%
\endbibitem

\bibitem{volz}
%
\begin{barticle}[author]
\bauthor{\bsnm{Volz},~\bfnm{E.}\binits{E.}}
(\byear{2008}).
\btitle{{S}{I}{R} dynamics in random networks with heterogeneous connectivity}.
\bjournal{J. Math. Biol.}
\bvolume{56}
\bpages{293--310}.
\end{barticle}
%
\MR{2358436}
\endbibitem

\bibitem{Whitt80}
%
\begin{barticle}[author]
\bauthor{\bsnm{Whitt},~\bfnm{W.}\binits{W.}}
(\byear{1985}).
\btitle{Blocking when service is required from several facilities
simultaneously}.
\bjournal{AT\&T Tech. J.}
\bvolume{64}
\bpages{1807--1856}.
\end{barticle}
%
\MR{0812939}
\endbibitem

\end{thebibliography}
\end{document}